\documentclass[11pt]{article}

\usepackage {amsfonts,amsthm}
\usepackage[english]{babel}

\theoremstyle{plain}
\newtheorem {Th}{Theorem}
\newtheorem*{Th0}{Theorem SB}
\newtheorem {Cor}{Corollary}
\newtheorem {Lem}{Lemma}
\newtheorem {Pro}{Proposition}
\newtheorem {Def}{Definition}
\newcommand {\bpr} {\noindent{\bf Proof.} }

\def\AE{{\cal AE}}
\def\NO{{\cal N}_1}
\def\card{{\rm{card}}}
\def\dist{{\rm{dist}}}
\def\diam{{\rm{diam}}}

\def\const{{\rm{const}}}

\def\C{{\mathbb{C}}}
\def\R{{\mathbb{R}}}
\def\Z{{\mathbb{Z}}}
\def\N{{\mathbb{N}}}
\def\Im{{\rm{Im\,}}}
\def\Re{{\rm{Re\,}}}
\def\s{\sigma}
\def\t{\tilde }
\def\wt{\widetilde }
\def\d{\delta}

\def\f{\varphi}
\def\a{\alpha}
\def\b{\beta}

\def\th{\theta}
\def\e{{\varepsilon}}
\def\l{\lambda}
\def\g{\gamma}
\def\supp{{\rm{supp}}}
\def\es{\emptyset}
\def\bs{~\hfill\rule{7pt}{7pt}}
\def\Dab{D=((a_n),(b_n))}

\def\be{\begin{equation}}
\def\ee{\end{equation}}

\begin{document}
\title{
    \bf Sunyer-i-Balaguer's Almost Elliptic Functions and Yosida's Normal Functions}

\author{Favorov S.Ju.}

\date{}

\maketitle

\begin{abstract}
We study the properties of two classes of meromorphic functions in
the complex plane. The first one is the class of almost elliptic
functions in the sense of Sunyer-i-Balaguer.  This is the class of
meromorphic functions $f$ such that the family
$\{f(z+h)\}_{h\in\C}$ is normal with respect to the uniform
convergence {\it in the whole complex plane}.  Given two sequences
of complex numbers, we provide sufficient conditions for them to
be zeros and poles of some almost elliptic function. These
conditions enable one to give (for the first time) explicit
non-trivial examples of almost elliptic functions.

 The second class was introduced by K.Yosida, who called it a class of normal
 functions of the first category. This is the class of meromorphic functions
$f$ such that the family $\{f(z+h)\}_{h\in\C}$ is normal with
respect to the uniform convergence {\it on compacta in the complex
plane} and no limit point of the family is a constant function. We
give necessary and sufficient conditions for two sequences of
complex numbers to be zeros and poles of some normal function of
the first category and obtain a parametric representation for this
class in terms of zeros and poles.
\end{abstract}



\bigskip

According to M.Bessonoff \cite{Be} and M.Favard \cite{F}, a
meromorphic function $f$ in the complex plain $\C$ is called
almost elliptic if the following condition holds: for each $\e>0$
and $\d>0$, there exists $L<\infty$ such that every real or pure
imaginary interval of the length $L$ contains points $\tau$ such
that the inequality
 $$
|f(z+\tau)-f(z)|<\e
 $$
holds for all points $z\in\C$ whose distance to the set of poles of
$f$ is larger than $\d$.

A subclass of such functions with a uniformly bounded number of
poles in all discs of radius 1 was investigated by H.Yoshida
\cite{Y2}.

 Another definition of an almost elliptic function was suggested by
F.Sunyer-i-Balaguer.
\begin{Def}[\cite{S}]\label{S}
A meromorphic function $f\not\equiv\const$ in the complex plain
$\C$ is called almost elliptic (we will say $f\in\AE$, or $f$ is
an $\AE$-function), if for every $\e>0$ there exists $L<\infty$
such that every real or pure imaginary interval of the length $L$
contains at most one $\e$-almost period of $f$, i.e., a point
$\tau$ with the property
 $$
\rho_S(f(z+\tau),f(z))<\e \quad \hbox{for all}\quad z\in\C,
 $$
where $\rho_S$ is the spherical metric in $\C$.
\end{Def}

In his seminal paper \cite{S} Sunyer-i-Balaguer investigated the
class $\AE$. In particular, he proved that $a$-points of every
$\AE$-function have  a uniform in a certain sense distribution.
Moreover, they are equidistributed for all $a\in\C\cup\{\infty\}$.
He also simplified conditions describing location of $\e$-almost
periods and proved Bochner's criterion

\medskip

\begin{Th0} The following conditions are equivalent

$\a$) $f\in\AE$,

$\b$) for each $\e>0$ there exists $L<\infty$ such that every disc
$\{z\in\C:\,|z-c|<L\}$ contains an $\e$-almost period $\tau$,

$\g$) for each sequence $(h_n)\subset\C$ there exists a
subsequence $(h_{n'})$ such that
$\rho_S(f(z+h_{n'}),f(z+h_{m'}))\to 0$ as $n',m'\to\infty$,
uniformly in $\C$.\footnote{See also Proposition \ref{1} below.}
\end{Th0}

Next, following to  K.Yosida \cite{Y1} (see also \cite{Y2},
\cite{M}, \cite{E}), we introduce the  definition of a normal
function.

\begin{Def}
A meromorphic function $f$ is called normal if  for each sequence
$(h_n)\subset\C$ there exists a subsequence $(h_{n'})$ such that
$\rho_S(f(z+h_{n'}),f(z+h_{m'}))\to 0$ as $n',m'\to\infty$
uniformly on compacta in $\C$.  A normal meromorphic function $f$
is of the first category ($f\in\NO$, or $f$ is an $\NO$-function)
if the family $\{f(z+h)\}_{h\in\C}$ has no constant functions as
limit points.
\end{Def}

A typical example of $\NO$-function is $1/z+E(z)$, where $E$ is an
arbitrary elliptic function.  The class of normal meromorphic
functions is well known and has proved to be useful in complex
analysis, and in particular, in the Nevanlinna theory (see
\cite{Y1}, \cite{Y2}, \cite{E}, \cite{LP}, \cite{M}, \cite{Z}).

However, until now, no description of zeros and poles of either
normal  or $\NO$-functions was  known.

The class $\AE$ formes a natural subclass of the class of normal
functions in $\C$, moreover, $\AE\subset\NO$. Yet, there were no
examples of $\AE$-functions except for usual elliptic ones. Note
that an appropriate description of the class of normal functions
in $\C^*=\C\setminus\{0\}$ (i.e., meromorphic functions $f$ such
that the family $\{f(\l z):\,\l\in\C^*\}$ is normal in $\C^*$) was
obtained by A.Ostrowski \cite{O} (see also \cite{Mo} and
\cite{E}).

 This paper is organized as follows.

In \S 1, we prove main properties of $\NO$-functions and
$\AE$-functions. Some of these properties are new, others were
described in \cite{E}, \cite{S}, \cite{Y2}. For reader's
convenience, below we will give complete proofs of these results.

In \S 2, we  obtain necessary and sufficient conditions for a pair
of discrete sets with multiplicities (in our definition, a
divisor) to be the zero set and the pole set, respectively, of an
$\NO$-function. We also show that any $\NO$-function is a
conditionally convergent meromorphic Weierstrass product of genus
1.

In \S 3, we introduce a general notion of an almost periodic
mapping. Following \cite{R} and \cite{FRR}, we introduce the
notion of an almost periodic divisor. We prove that each
$\AE$-function has an almost periodic divisor. Next, under an
additional assumption on the regularity of indexing  of the
divisor, we investigate properties of $\arg f$.

In \S 4, we find necessary and sufficient conditions for an almost
periodic divisor with a regular indexing to be a divisor of an
$\AE$-function $f$. Moreover, if this is the case, then $f$
 is a conditionally convergent meromorphic Weierstrass product of genus 0.

In \S 5, we construct two examples of $\AE$-functions, which are
not elliptic.

  {\bf Acknowledgements.} I am very grateful to
professor A.Eremenko from whom I learned about these problems.

\bigskip
\begin{center}
{\bf \S 1. Main properties of $\NO$-functions and $\AE$-functions}
\end{center}
\bigskip

Let $\Z$ be the set of all integers. A mapping $D:\,\C\to\Z$ is a
{\it divisor} in the complex plane $\C$ if $\supp D$ is a discrete
set; the divisor is {\it positive} if $D(\C)\subset\Z^+\cup\{0\}$.

Let $(a_n)$ be a sequence of points from the set
$D^{-1}\{k\in\Z:\,k>0\}$, in which every point $a$ occurs $D(a)$
times, and $(b_n)$ be a sequence of points from the set
$D^{-1}\{k\in\Z:\,k<0\}$, in which every point $b$ occurs $-D(b)$
times. We will write a divisor $D$ as $((a_n),(b_n))$. Also,
$((a_n),\es)$ means a positive divisor. Furthermore, $\Dab$ is the
divisor of a meromorphic function $f$ if $\{a_n\}$ is the zero set
of $f$, $\{b_n\}$ is the pole set of $f$, $D(a_n)$ are
multiplicities of the zeros $a_n$, and $-D(b_n)$ are
multiplicities of the poles $b_n$.  By $\mu_D$ we denote the
discrete signed measure with masses $D(a_n)$ at the points $a_n$
and negative masses $D(b_n)$ at the points $b_n$. Let $m_2(E)$ be
the two-dimensional Lebesgue measure of a set $E\subset\C$,
$m_1(E)$ be the Hausdorff 1-dimensional measure of $E$. If $E$ is
a rectifiable curve, then $m_1(E)$ is just the length of $E$. By
$B(c,R)$ denote the disc $\{z\in\C:\,|z-c|<R\}$, by $\card A$
denote the number of points in a discrete set $A$, and by $C$
without indices denote any constant depending only on the divisor
$D$ or the function $f$.

 We begin with some simple properties of
 $\NO$-functions and $\AE$-functions (for $\AE$-functions see
\cite{S}).

It is easy to check that if $f$ is $\NO$-function  and for a
sequence $\{h_n\}$
 $$
\rho_S(f(z+h_n),f(z+h_m))\to 0\quad\hbox{as}\quad n,m\to\infty
 $$
uniformly on compacta in $\C$, then there exists a function
$g\in\NO$ such that
 $$
 \rho_S(f(z+h_n),g(z))\to 0\quad\hbox{as}\quad n\to\infty
 $$
 uniformly on compacta in $\C$. Moreover, the classes $\NO$ and $\AE$ are
 closed with respect to the uniform convergence in $\C$.
 It is also clear that
every $\NO$-function is a uniformly continuous map of $\C$ into
the Riemann sphere. Next, a composition of a rational function
with an $\NO$-function ($\AE$-function) is the $\NO$-function
($\AE$-function).

\begin{Th}[\cite{Y2}]\label{8}
The divisor of any $\NO$-function $f$ satisfies the separation
condition
\begin{equation}\label{sep}
 \inf_{n,k}|a_n-b_k|=\d_0>0.
  \end{equation}
\end{Th}
\bpr If $f(a_n)=0$, then $|f(z)|<1$ for $z\in B(a_n,\d_0)$, where
$\d_0$ does not depend on $n$.\bs

\begin{Th}[\cite{Y2}]\label{10}
Suppose $\Dab$ is the divisor of an $\NO$-function $f$, and $\d>0$
is an arbitrary real number. Then there exists $C_\d<\infty$ such
that
 \begin{equation}\label{p}
 |f(z)|\le C_\d\quad\forall z\not\in\bigcup_nB(b_n,\d),
 \end{equation}
 \begin{equation}\label{z}
 |f(z)|\ge 1/C_\d\quad\forall z\not\in\bigcup_nB(a_n,\d),
 \end{equation}
and
 \begin{equation}\label{der}
 |f'(z)/f(z)|\le C'_\d \quad\forall
 z\not\in\bigcup_nB(b_n,\d)\cup\bigcup_nB(a_n,\d).
 \end{equation}
 \end{Th}

\bpr Suppose (\ref{z}) is false. Take a sequence $(w_k)$ such that
$w_k\not\in\cup_n B(a_n,\d)$ for all $k$ and $f(w_k)\to 0$ as
$k\to\infty$. Then there is a subsequence $w_{k'}$ such that
$\rho_S(f(z+w_{k'}),\,g(z))\to0$ uniformly on compacta in $\C$. In
particular, $g(0)=0$. Now Hurwitz' Theorem leads to a
contradiction. In the same way we obtain (\ref{p}). Using
(\ref{z}) and (\ref{p}) for $\d/2$, we get (\ref{der}). \bs

\begin{Th}[\cite{Y2}]\label{9}
Let $\Dab$ be the divisor of an $\NO$-function $f$, and $|\mu_D|$
be the variation of the measure $\mu_D$. Then
 \begin{equation}\label{d}
|\mu_D|(\overline{B(c,1)})=\card\{n:\,a_n\in
\overline{B(c,1)}\}+\card\{n:\,b_n\in
\overline{B(c,1)}\}<C_0\quad\forall
 c\in\C.
 \end{equation}
\end{Th}

\bpr If  $\card\{n:\,a_n\in B(c_k,\d_0/3)\}\to\infty$ as
$k\to\infty$,  then for a subsequence $(c_{k'})$ we get
$\rho(f(z+c_{k'}),\,g(z))\to 0$ uniformly on compacta in $\C$,
where $g$ is a nonzero meromorphic function. In view of
(\ref{sep}) and (\ref{p}), we get $|f(z+c_k)|<C$ for $|z|<\d_0/3$,
hence, $|f(z+c_{k'})-g(z)|\to 0$ in this disc. Then Hurwitz'
theorem leads to a contradiction. Therefore, $\card\{n:\,a_n\in
B(c,\d_0/3)\}<C$ for all $c\in\C$. The same is valid for poles of
$f$. Thus, we obtain (\ref{d}).\bs

\begin{Cor}\label{9a}
Suppose $D$ is a divisor with property (\ref{d}); then
 \begin{equation}\label{d1}
|\mu_D|(E)\le C_0S(E),\quad\hbox{for every set}\quad E\subset\C,
 \end{equation}
where $S(E)$ is the number of closed discs of radius 1 covering
$E$,
\begin{equation}\label{m1}
  |\mu_D|(B(c,r))<Cr^2,\quad  \forall c\in\C,\quad r>1,
  \end{equation}
and
  \begin{equation}\label{m2}
   \int_{r<|w-c|}|w-c|^{-k}d|\mu_D|(w)<Cr^{-(k-2)},
   \quad\forall c\in\C,\quad r>1,\quad k>2.
 \end{equation}
 \end{Cor}
\bpr Clearly, (\ref{m1}) follows from (\ref{d1}), and implies
(\ref{m2}). \bs

\medskip
Theorem \ref{10}, Proposition \ref{i}, and Liouville's Theorem
yield

 \begin{Cor}\label{ent}
Every entire $\NO$-function is a constant.
 \end{Cor}

 \begin{Cor}\label{11a}
Suppose $f_1,\,f_2$ are $\NO$-functions with the same divisor;
then $f_1=Kf_2$ with a constant $K\in\C$.
 \end{Cor}

Note that Theorem \ref{10}, Theorem \ref{9}, and Corollaries
\ref{9a} -- \ref{11a} for $\AE$-functions were proved in \cite{S}.

 \medskip
The following simple proposition are needed for the sequel.
\begin{Pro}\label{i}
Suppose $D$ is a divisor with property (\ref{d}) and
$\d<1/(2C_0)$. Then the diameter of every connected component $A$
of the set $A(\d)=(\cup_n B(a_n,\d))\cup(\cup_n B(b_n,\d))$ is at
most $2\d C_0$.
\end{Pro}
\bpr It follows from (\ref{d1}) that any circle $\{z:\,|z-c|=1\}$,
$c\in A$, has no common point with $A$. Therefore the disc
$B(c,\,1)$ contains $A$. The statement now follows from
(\ref{d}).\bs

\begin{Th}[\cite{Y1}]\label{7a}
For any $f\in\NO$ with a divisor $\Dab$ there is $R_0$ such that
\begin{equation}\label{N1}
  B(w,R_0)\cap\{a_n\}_{n\in\N}\neq\emptyset,\quad
  B(w,R_0)\cap\{b_n\}_{n\in\N}\neq\emptyset,
  \quad \forall\, w\in\C.
\end{equation}
\end{Th}
\bpr   Let $B(w_k,R_k)$ be discs without poles of $f$ such that
$R_k\to\infty$. Using (\ref{p}) with $\d=1$, we see that
$|f(w_k+z)|<C$ for $|z|<R_k-1,\,k=1,\dots$. Furthermore,
$f(z+w_{k'}) \to g(z)$ uniformly on compacta in $\C$ for a
subsequence $(w_{k'})\subset(w_{k'})$. Hence we get
$g(z)\equiv\const$, which is impossible. Since $1/f\in\NO$, we
obtain the similar result for zeros of $f$. \bs

\begin{Th}[\cite{Y2}]\label{Y}
A meromorphic function $f$ is an $\NO$-function if and only if the
following conditions are fulfilled:

a) for any $\d>0$
\begin{equation}\label{pz}
 1/C_\d\le|f(z)|\le C_\d,\quad\forall
 z\not\in A(\d)=\bigcup_nB(a_n,\d)\cup\bigcup_nB(b_n,\d),
 \end{equation}

b) zeros and poles of $f$ satisfies  (\ref{sep}), (\ref{d}), and
(\ref{N1}).
\end{Th}

\bpr For any $\NO$-function $f$ properties a) and b) follows from
Theorems \ref{8}, \ref{9}, \ref{10}, and \ref{7a}.

Let a meromorphic function $f$ satisfy a) and b). If $\d_0$ is the
 same as in (\ref{sep}), then any connected component $A$ of the
 set $A(\d_0/3)$ can not  contain zeros and poles simultaneously.
 If it does  not contain poles of $f$, then inequality (\ref{pz})
  and the Maximum Modulus Principle yield
 $|f(z)|<C$ for $z\in A$. Similarly, if $A$
does not contain zeros of  $f$, then $|f(z)|>1/C$ for $z\in A$.
Take a disc $B(z_0,\d_0/3)$. We see that for each sequence $(w_k)$
 there is a subsequence $(w_{k'})$ such that at points of the disc
either $|f(z+w_{k'})|<C$ for all $k'$, or $|1/f(z+w_{k'})|<C$ for
all $k'$. Hence, in both cases there is a subsequence
$(w_{k''})\subset(w_{k'})$ and a meromorphic function $g(z)$ in
the disc $B(z_0,\d_0/3)$ such that uniformly in $z\in
B(z_0,\d_0/4)$
\begin{equation}\label{conv}
\rho_S(f(z+w_{k''}),g(z))\to 0\quad\hbox{as}\quad k''\to\infty.
\end{equation}
 Let $\{B(z_n,\d_0/4)\}$ be a denumerable covering of $\C$.
 Using the diagonal procedure, we obtain the subsequence
 $(w_{k'''})\subset(w_{k''})$ such that
(\ref{conv}) holds uniformly on compacta in $\C$.

In view of (\ref{N1}), every function $f(z+w_{k''})$ has a pole
and a zero in the disc $B(0,R_0)$. Hence, $g(z)\not\equiv\const$
and $f\in\NO$. \bs

\medskip
The same argument proves the following theorem:
\begin{Th}\label{Y1}
If a meromorphic function $f$ satisfies (\ref{pz}) for any $\d>0$,
and its divisor satisfies (\ref{sep}) and (\ref{d}), then $f$ is a
normal function.
\end{Th}

 \begin{Th}\label{11}
The product of $f_1,\,f_2\in\NO$is a $\NO$-function if and only if
the divisor of the function $f_1f_2$ satisfies (\ref{sep}) and
(\ref{N1}). A similar assertion holds for the class $\AE$.
 \end{Th}

\bpr Necessity follows from Theorems \ref{8} and \ref{7a}. To
prove sufficiency, take a sequence $(w_k)\subset\C$. Then there
exist functions $g_1,\,g_2\in\NO$ and a subsequence $(w_{k'})$
such that
 $$
\rho_S (f_1(z+w_{k'}),\, g_1(z))\to 0,\quad
 \rho_S (f_2(z+w_{k'}), \,g_2(z))\to 0\quad \hbox{as}\quad
 k'\to\infty,
 $$
 uniformly on compacta in $\C$.
 Let $U$ be the union of the discs $B(b_j,\d)$ over all
 poles $b_j$ of $g_1$ and $g_2$. Using Theorem
 \ref{10}, we obtain that uniformly on compacta in $\C\setminus U$
 $$
f_1(z+w_{k'})-g_1(z)\to 0,\quad f_2(z+w_{k'})-g_2(z)\to 0,
 $$
 \begin{equation}\label{pr1}
 (f_1f_2)(z+w_{k'})-(g_1g_2)(z)\to 0\quad\hbox{as}\quad k'\to\infty.
  \end{equation}
Suppose that the distances between zeros and poles of the function
$f_1f_2$ are at least $\e$. Taking into account Proposition
\ref{i}, we see that for sufficiently small $\d$, the diameter of
any connected component $A$ of the set $U$ is less than $\e$.
Hence, $A$ does not contain simultaneously poles and zeros of the
function $(f_1f_2)(z+w_{k'})$. If $A$ does not contain poles of
$(f_1f_2)(z+w_{k''})$ for a subsequence
$(w_{k''})\subset(w_{k'})$, then the Maximum Modulus Principle and
(\ref{pr1}) imply the convergence of the functions
$(f_1f_2)(z+w_{k''})$ to $(g_1g_2)(z)$ uniformly in $z\in A$. If
$A$ does not contain zeros of $(f_1f_2)(z+w_{k''})$, then the same
argument shows that
 $$
 1/(f_1f_2)(z+w_{k''})-1/(g_1g_2)(z)\to 0\quad\hbox{as}\quad k'\to\infty
 $$
uniformly in $z\in A$. Consequently,
$$
\rho_S (f_1f_2(z+w_{k''}),\,g_1g_2(z))\to 0\quad \hbox{as}\quad
 k'\to\infty
 $$
uniformly on compacta in $\C$. In view of (\ref{N1}), we get
$g_1g_2\not\equiv\const$ and $f_1f_2\in\NO$. Likewise, if
$f_1,\,f_2\in\AE$ and the divisor of $f_1f_2$ satisfies
(\ref{sep}) and (\ref{N1}), then $f_1f_2\in\AE$. \bs

\begin{Th}[\cite{E}]\label{20}
For any $\NO$-function $f$ there is a constant $C<\infty$ such
that
 \begin{equation}\label{I}
   \int_0^{2\pi}|\log|f(c+re^{i\th})||d\th<C,
 \end{equation}
 for all $c\in\C$ and $r>1$.
\end{Th}

\bpr Let $\Dab$ be the divisor of $f$, and $A_1,\dots,A_N$ the
connected components of the set
$A(\d)=\bigcup_nB(a_n,\d)\cup\bigcup_nB(b_n,\d)$, which have
nonempty intersection with the circle $z=c+re^{i\th}$. Put
$\d=\min\{\d_0/(2C_0),1/(2C_0)\}$. By Proposition \ref{i}, it
follows that $\diam A_j<\min\{1,\d_0\}$ for every $j=1,\dots,N$.
From (\ref{sep}) it follows that each $A_j$ does not contain poles
and zeros of $f$ simultaneously. From (\ref{d1}) it also follows
that the number $k_j$ of zeros or poles in the set $A_j$ does not
exceed $C_0$. Since the annulus $r-1<|w|<r+1$ is covered by $6\pi
r$ discs of radius $1$, we have
 \begin{equation}\label{k}
   \sum_1^N k_j<6\pi C_0r.
 \end{equation}
 Furthermore, using (\ref{z}) and (\ref{p}), we get
  \begin{equation}\label{II}
 \int_0^{2\pi}|\log|f(c+re^{i\th})||d\th \le
 \sum_1^N \int_{E_j}|\log|f(c+re^{i\theta})||d\theta
 +2\pi\log C,
  \end{equation}
 where $E_j=\{\th:\,c+re^{i\th}\in A_j\}$. Suppose a
 component $A_j$ contains zeros of $f$.  Take
 $P_j(z)=\d^{-k_j}\prod_{a_n\in A_j}(z-a_n)$. Since
 $\diam A_j< \d_0$, we get
 $$
 1\le |P_j(z)|\le (\diam A_j/\d)^{k_j}<(2C_0)^{C_0}\qquad\forall z\in\partial A_j.
 $$
Therefore, $|\log|f(z)/P_j(z)||\le C$ for $z\in\partial A_j$. The
Maximum Modulus Principle yields the same inequality for all $z\in
A_j$. Hence,
\begin{equation}\label{I1}
\int_{E_j}|\log|f(c+re^{i\theta})||d\theta \le \int_{E_j}
|\log|P_j(c+re^{i\theta})||d\theta +Cm_1(E_j).
\end{equation}
Next, we have
\begin{equation}\label{I2}
\int_{E_j} \log^+|P_j(c+re^{i\theta})|d\theta \le
m_1(E_j)C_0\log(2C_0),
\end{equation}
and
\begin{equation}\label{I3}
\int_{E_j} \log^+\left|\frac{1}{P_j(c+re^{i\theta})}\right|d\theta
\le \sum_{a_n\in A_j}\int_{\th:|re^{i\th}-a_n+c|<\d}
\log^+\left|\frac{\d}{re^{i\th}+c-a_n}\right| d\th.
\end{equation}
Note that for every $te^{i\f}\in\C,\,t\ge0$, we have
 $$
\int_{\th:|re^{i\th}-te^{i\f}|<\d}
\log^+\left|\frac{\d}{re^{i\th}-te^{i\f}}\right| d\th <
\int_{|\th-\f|<\pi\d/2r}
\log^+\left|\frac{\d/r}{(\th-\f)2/\pi}\right|d\th <\d\pi/2r.
 $$
Combining (\ref{I1}), (\ref{I2}), and (\ref{I3}), we obtain
 $$
\int_{E_j} |\log^+|f(c+re^{i\theta})||d\theta < C(m_1(E_j)+k_j/r).
 $$
A similar bound is valid for $A_j$ containing poles of $f$.
Therefore, (\ref{II}) and (\ref{k}) imply (\ref{I}). \bs

\begin{Th}\label{c3}
 Suppose $f\in\NO$ has the divisor $\Dab$; then
  $$
\left|\int_1^r{\card\{n:\,|a_n|<t\}-\card\{n:\,|b_n|<t\}\over t}
\,dt\right|<C.
 $$
The constant $C$ is the same for all $r>1$ and all shifts of $f$.
\end{Th}
\bpr Using Jensen's formula for the function $f$, we get
 $$
\int_1^r{\mu_D(B(c,t))\over t}\,dt={1\over 2\pi}
\int_0^{2\pi}\log|f(c+re^{i\th})|d\th-{1\over
2\pi}\int_0^{2\pi}\log|f(c+e^{i\th})|d\th.
  $$
The result now follows from (\ref{I}).
 \bs
\medskip

Theorem \ref{c3} has an analogue for algebraic sums of zeros and
poles.
 \begin{Th}\label{6}
Suppose $f\in\NO$ has the divisor $\Dab$; then
 \begin{equation}\label{Int}
\left|{1\over r}\int_1^r
\left(\sum_{n:\,|a_n|<t}a_n-\sum_{n:\,|b_n|<t}b_n\right) {dt\over
t}\right|<C\qquad\forall\,r>1.
\end{equation}
The constant $C$ is the same for all shifts of $f$.
 \end{Th}
\bpr Put $w=u+iv$. Consider the distribution
$\triangle(u\log|f(w)|)$. Since the function $(\log|f(w)|)'_u$ is
locally integrable over $\C$, we see that  the function
$u\log|f(w)|$ up to a harmonic function is a logarithmic potential
of the measure $u\cdot2\pi\mu_D(w)+2(\log|f(w)|)'_u m_2(w)$ in the
disc $B(0,R)$. Applying Jensen--Privalov's formula for the annulus
$1\le|w|\le r$, we get
 $$
{1\over 2\pi}
\int_0^{2\pi}(r\cos\th)\log|f(re^{i\th})|d\th-{1\over
2\pi}\int_0^{2\pi}\cos\th\log|f(e^{i\th})|d\th
 $$
 $$
 =\int_1^r\left[\int_{B(0,t)}u~d\mu_D(w)+\int_{B(0,t)}2(\log|f(w)|)'_u~dm_2\right]
 {dt\over t}.
 $$
 Note that
 $$
 \int_{B(0,t)}(\log|f(w)|)'_u
 ~dm_2(w)=\int^t_{-t}(\log|f(\sqrt{t^2-v^2}+iv)|-\log|f(-\sqrt{t^2-v^2}+iv)|)~dv
 $$
 $$
 =\int^{2\pi}_0\log|f(t\cos\th+it\sin\th)|t\cos\th~d\th.
 $$
 Taking into account (\ref{I}), we obtain
 $$
\left|{1\over r}
\int_1^r\left[\int_{B(0,t)}u~d\mu_D(w)\right]{dt\over t}\right|<C.
 $$
 The same bound is valid for the measures
$v\cdot\mu_D(w)$ and  $(u+iv)\cdot\mu_D(w)$. Hence, we obtain
(\ref{Int}). Obviously,  the constant $C$ does not depend on
shifts of $f$. \bs
\medskip

 \begin{Th}\label{tr}
Let $\Dab$ be the divisor of an $\NO$-function $f$. Then for any
simply connected bounded domain $E\subset\C$ such that
$m_1(\partial E)<\infty$, and any holomorphic function $g$ in $E$
we have
\begin{equation}\label{es}
   \left|\int_E g(w)\;d\mu_D(w)\right|\le C\sup_E|g(w)|(m_1(\partial
   E)+1).
 \end{equation}
In particular, putting $g(w)\equiv 1$, we obtain
 \begin{equation}\label{es1}
   |\mu_D(E)|\le C(m_1(\partial E)+1).
 \end{equation}
Next, for any $k\in\Z$ and $1<r<R<\infty$
 \begin{equation}\label{es2}
   \left|\int_{r<|w|\le R} w^k\;d\mu_D(w)\right| <
   C(r^{k+1}+R^{k+1}).
 \end{equation}
The constants in (\ref{es}), (\ref{es1}), and (\ref{es2}) do not
depend on $g$ and shifts of $f$.
\end{Th}

We need the following simple lemma.
 \begin{Lem}\label{new}
   Suppose a compact set $F\subset\C$  consists of $M$ connected
components and $m_1(F)<\infty$. Then there are at most
$M(6m_1(F)+1)$ discs $B(z_k,1)$ such that $\cup_kB(z_k,1)\supset
F$.
 \end{Lem}

\bpr Let $F$ be a connected set. If $F\subset B(z_0,\,1)$ for some
$z_0\in F$, there is nothing to prove. Otherwise, by Besicovitch's
covering principle (\cite{Besi}, also see \cite{La}, Lemma 3.2),
there are discs $B(z_k,1),\,z_k\in F$, $k=1,\dots,S$, such that
$\cup_kB(z_k,1)\supset F$ and every point of $\C$ belongs to at
most $6$ discs. Clearly, for each disc we have $m_1(F\cap
B(z_k,1))\ge 1$. Therefore\footnote{$m_1(E)$ is an additive
function on Borel sets, see \cite{DF}, p.191}, $S\le\sum_k
m_1(F\cap B(z_k,1))\le 6m_1(F)$.

For the general case, one should apply this result to each
connected component of $F$.\bs
\medskip

\noindent{\bf Proof of the Theorem.} Using Lemma \ref{new}, take a
covering of $\partial E$ by discs $B(z_k,1),\,k=1,\dots,S,\,S\le
6(m_1(\partial E)+1)$. Put $E_1=E\setminus\overline{\cup_k^S
B(z_k,1)}$. We have
 $$
   m_1(\partial E_1) \le2\pi S \le12\pi(m_1(\partial E)+1).
 $$
Next, using (\ref{d1}), we get
\begin{equation}\label{B2}
   |\mu_D|(E\setminus E_1)\le6C_0(m_1(\partial E)+1).
 \end{equation}
Take $\d=1/(2C_0)$. Let $A_1,\dots,A_N$ be all connected
components of the set $A(\d)=\cup_nB(a_n,\d)\cup\cup_nB(b_n,\d)$
with nonempty intersection with $\partial E_1$. Clearly, the
diameter of every connected component is less than $1$. Therefore,
$\overline{\cup_j^N A_j}\subset\cup_1^S B(z_k,2)$.  Using
(\ref{m1}), we obtain
 \begin{equation}\label{B3}
|\mu_D|(\overline{\cup_j^N A_j})\le|\mu_D|(\cup_1^S B(z_k,2))\le
C(m_1(\partial E)+1).
 \end{equation}
 Moreover,
  $$
 m_1(\cup_j^N\partial A_j)\le 2\pi\d|\mu_D|(\overline{\cup_j^N A_j})
 \le C(m_1(\partial E)+1).
 $$
Put $E_2=E_1\setminus\overline{\cup_j^N A_j}$. Obviously, $E_2$ is
a finite union of domains in $\C$, and $\partial E_2$ is a finite
union of circular arcs such that
 $$
   m_1(\partial E_2)\le m_1(\partial E_1)+ m_1(\cup_j^N\partial A_j)\le C(m_1(\partial E)+1).
 $$
In addition, $\partial E_2\cap A(\d)=\es$. Using (\ref{der}), we
obtain that right-hand side of the equality
 \begin{equation}\label{m}
   \int_{E_2} g(w)\;d\mu_D(w)={1\over 2\pi i}\int_{\partial E_2}
{g(w)f'(w)\over f(w)}dw
 \end{equation}
 does not exceed $C\max_{\partial E_2}|g(w)|(m_1(\partial E)+1)$.
This bound together with (\ref{B2}) and (\ref{B3}) proves
(\ref{es}).

To prove (\ref{es2}), we apply the above argument to the
doubly-connected domain $E=\{w:\,r<|w|\le R\}$ with $g(w)=w^k$.
The integral in the right-hand side of (\ref{m}) over $\partial
E_2\cap\{w:\,|w|<r+3\}$ is bounded by $Cr^{k+1}$ and the integral
over $\partial E_2\cap\{w:\,|w|>R-3\}$ is bounded by $CR^{k+1}$.
\bs
\medskip

Putting $k = -2$ and $k = -1$ in (\ref{es2}), we obtain
\begin{Cor}\label{c1}
For  every $\NO$-function with the divisor $\Dab$ there is a
finite limit
 $$
\lim_{r\to\infty}\sum_{1\le
|a_n|<r}1/a_n^2-\sum_{1\le|b_n|<r}1/b_n^2.
 $$
Next, for all $r>1$
  $$
\left|\sum_{1\le |a_n|<r}1/a_n-\sum_{1\le|b_n|<r}1/b_n\right|<C.
  $$
The constant $C$ does not depend on shifts of $D$.
\end{Cor}

{\bf Remark}. In fact, we have proved all theorems of this section
(with the exception of Theorems \ref{7a} and \ref{Y}) for normal
functions $f$ such that $g(z)\equiv 0$ and $g(z)\equiv\infty$ are
not limiting functions for the family of shifts of $f$.

\bigskip
\begin{center}
{\bf \S 2. Representation for $\NO$-functions and description of
zero sets and pole sets}
\end{center}
\bigskip
The main result of this section is
\begin{Th}\label{Rep}
 A divisor $\Dab,\,
a_n\neq 0,\,b_n\neq 0$ for all $n$ is the divisor of $f\in\NO$, if
and only if  the following conditions are fulfilled:
\medskip

a) $\inf_{n,k}|a_n-b_k|>0$,
\medskip

b) $\card\{n:\,|a_n|\le1\}+\card\{n:\,|b_n|\le1\}<C$, uniformly
with respect to shifts of $f$,
\medskip

c) there exists a radius $R<\infty$ such that every disc $B(c,R)$,
$c\in\C$, intersects with $\{a_n\}$ and $\{b_n\}$, simultaneously,
\medskip

d) $$
\left|\int_1^r{\card\{n:\,|a_n|<t\}-\card\{n:\,|b_n|<t\}\over t}
\,dt\right|<C
 $$
 for all $r>1$ uniformly with respect to shifts of $f$,

e) $$
 \left|{1\over r}\int_1^r
\left(\sum_{n:\,|a_n|<t}a_n-\sum_{n:\,|b_n|<t}b_n\right) {dt\over
t}\right|<C,
  $$
for all $r>1$ uniformly with respect to shifts of $f$,

 f) there exists a finite limit
 $$
\lim_{r\to\infty}\sum_{1\le
|a_n|<r}1/a_n^2-\sum_{1\le|b_n|<r}1/b_n^2.
 $$

Moreover, each $\NO$-function with the divisor $\Dab$ up to a
constant factor has the form
\begin{equation}\label{rep}
  f(z)=e^{\a z}\lim_{r\to\infty}{\prod_{n:|a_n|<r}(1-z/a_n)e^{z/a_n}\over
 \prod_{n:|b_n|<r}(1-z/b_n)e^{z/b_n}}.
\end{equation}
Here the limit exists uniformly on compacta in $\C$ and
\begin{equation}\label{g}
 \a=\lim_{r\to\infty}
\sum_{n:|b_n|<r}(1/b_n-\overline{b_n}/r^2)-\sum_{n:|a_n|<r}(1/a_n-\overline{a_n}/r^2).
\end{equation}
\end{Th}

\bpr For a function $f\in\NO$,  conditions a), b), c), d), e), and
f) follow from Theorems \ref{8}, \ref{9}, \ref{7a}, \ref{c3},
\ref{6}, and Corollary \ref{c1}, respectively.

Let us prove sufficiency of these conditions. Put
 $$
 \a(r)=\int_0^r(1/|w|^2-1/r^2)\overline{w}\,d\mu_D(w),\quad
 \b(t)=\int_{B(0,t)}\overline{w}\,d\mu_D(w),\quad
 \g(t)=\int_0^t(\b(s)/s)\,ds.
 $$
Note that $\supp\mu_D\cap B(0,t)=\es$ for sufficiently small $t$.
Integrating by parts, we get
 $$
\a(R)-\a(r)=2\int_r^R{\b(t)\over t^3}\,dt=
 {2\g(R)\over R^2}
 -{2\g(r)\over r^2}+4\int_r^R{\g(t)\over t^3}\,dt.
 $$
 Using e), we get $|\g(t)|<Ct$.
Therefore, the limit $\a=-\lim_{r\to\infty}\a(r)$ in (\ref{g})
exists.

Furthermore, condition b) yields (\ref{d1}) and (\ref{m1}).
Conditions (\ref{m1}) and c) imply that both sequences $(a_n)$ and
$(b_n)$ have the genus $2$. Using f), we obtain that the function
$f$ in (\ref{rep}) is well defined. Let us show that it belongs to
the class $\NO$.

 It follows from (\ref{g}) that we can rewrite (\ref{rep}) in the form
 $$
  f(z)=\lim_{r\to\infty}f_r(z),
 $$
where
 $$
f_r(z)={\prod_{n:|a_n|<r}(1-z/a_n)e^{\bar a_nz/r^2}\over
 \prod_{n:|b_n|<r}(1-z/b_n)e^{\bar b_nz/r^2}}.
 $$
Integrating by parts, we get
 $$
 \log|f_r(z)|=\int\limits_{|w|<r}\log\left|{z-w\over r}\right|\;d\mu_D(w)
-\int\limits_{|w|<r}\log\left|{w\over
r}\right|\;d\mu_D(w)+\Re\left\{\int\limits_{|w|<r}{z\overline{w}\over
r^2}\;d\mu_D(w)\right\}
 $$
  \begin{equation}\label{b0}
 = \int\limits_0^r{\mu_D(B(0,t))-\mu_D(B(z,t))\over t}dt+
 \int\limits_{|w|<r, |w-z|\ge r}\log\left|{w-z\over
r}\right|\;d\mu_D(w)
 \end{equation}
 $$
  -\int\limits_{|w-z|<r, |w|\ge
r}\log\left|{w-z\over r}\right|\;d\mu_D(w)+
\Re\left\{\int\limits_{|w|<r}{z\overline{w}\over
r^2}\;d\mu_D(w)\right\}.
 $$
Let $s(w)$ be a nonnegative number such that $|w-s(w)z|=r$. If
$|w|<r$ and $|w-z|\ge r$, then $s(w)\in(0,\,1]$ and
 $$
\log\left|{w-z\over r}\right|=\Re{z(s(w)-1)\over
(w-s(w)z)}+O(|z|^2/r^2)
 =\Re{z\overline{w}(s(w)-1)\over r^2}+O(|z|^2/r^2)
\quad\hbox{as}\quad r\to\infty.
 $$
Since $\mu_D$ satisfies (\ref{d1}), we see that
 $$
|\mu_D|(\{|w|<r, |w-z|\ge r\})\le Cr(|z|+1).
 $$
 Hence,
 $$
\int_{|w|<r, |w-z|\ge r}\log\left|{w-z\over r}\right|\;d\mu_D(w)
=\Re\left\{{z\over r^2}\int_{|w|<r, |w-z|\ge
r}(s(w)-1)\overline{w}\,d\mu_D(w)\right\}+o(1)
 $$
as $r\to\infty$. An application of the same argument shows that
the sum of the last three integrals in (\ref{b0}) is equal to the
real part of the sum
 \begin{equation}\label{b}
{z\over r^2}\left[\int\limits_{|w-z|<r}\overline{w}\;d\mu_D(w)
+\int\limits_{|w|<r, |w-z|\ge
r}s(w)\overline{w}\;d\mu_D(w)-\int\limits_{|w|\ge r,
|w-z|<r}s(w)\overline{w}\;d\mu_D(w)\right]
 \end{equation}
up to the term $o(1)$ as $r\to\infty$. Note that every point $w\in
B(0,r)\setminus B(z,r)$ belongs to $B(sz,r)$ only for
$s\in[0,s(w))$ and every point $w\in B(z,r)\setminus B(0,r)$
belongs to $B(sz,r)$ only for $s\in(s(w),1]$. Therefore, (\ref{b})
is equal to the integral
 $$
{z\over r^2}\int_0^1 \int_{|w-sz|<r}\overline{w}\;d\mu_D(w)~ds.
 $$
Hence we obtain
\begin{equation}\label{log}
\log|f(z)|=
\lim_{r\to\infty}\int\limits_0^r{\mu_D(B(0,t))-\mu_D(B(z,t))\over
t}dt +\Re\left(\int\limits_0^1\int\limits_{|w-sz|<r}
{z\overline{w}\over r^2}~d\mu_D(w)\,ds\right).
\end{equation}

 Take $\d\in (0,\,1)$ such that $0\not\in A(\d)$. In view of b), the
integral
 $$
\int_0^1 {\mu_D(B(0,t))-\mu_D(B(z,t))\over t}dt
 $$
is uniformly bounded in $z\in\C\setminus A(\d)$. Also, by d), the
integral
 $$
\int\limits_1^r{\mu_D(B(0,t))-\mu_D(B(z,t))\over t}dt
 $$
is  uniformly bounded in $z\in\C$ and $r>1$ as well. Furthermore,
since bound e) does not depend on shifts of $\mu_D$, we get for
all $z\in\C$ and $r<R$,
 $$
\left|\int_r^R{1\over t}\int_{B(z,t)}\overline{w-z}~
d\mu_D(w)\,dt\right|<C(r+R).
 $$
In view of d), we get
 $$
\left|\int_r^R{1\over t}\int_{B(z,t)}d\mu_D(w)\,dt\right|<C.
 $$
Therefore,
 $$
\Re\left({z\over|z|}\int_r^R{1\over t}\int_{B(z,t)}\overline{w}~
d\mu_D(w)\,dt\right)\le \left|\int_r^R{1\over t}
\int_{B(z,t)}\overline{w}~d\mu_D(w)\,dt\right|<C(r+R+|z|).
 $$
Replace $z$ by $sz$, $R$ by $r+|z|$, and integrate over $s$ from
$0$ to $1$. We get
 $$
{1\over|z|} \int_r^{r+|z|}\Re\left({z\over
t}\int_0^1\int_{B(sz,t)}\overline{w}\,d\mu_D(w)\,ds\right)dt<C(r+|z|).
 $$
Therefore, for some $r'\in(r,\,r+|z|)$
 $$
\Re\left({z\over r'} \int_0^1\int_{B(sz,r')}\overline{w}\,
d\mu_D(w)\,ds\right)<C(r+|z|).
 $$
Hence for a sequence $r'\to\infty$
 $$
\Re\left({z\over (r')^2}\int\limits_0^1\int\limits_{|w-sz|<r'}
\overline{w}\;d\mu_D(w)\,ds\right)\le 2C.
 $$
 Similarly, for some sequence $r''\to\infty$
 $$ \int\limits_0^1\Re\left({z\over
(r'')^2}\int\limits_{|w-sz|<r''}
\overline{w}\;d\mu_D(w)\right)ds\ge -2C.
 $$
Taking into account (\ref{log}) and d), we get (\ref{pz}). Now, by
Theorem \ref{Y}, $f\in\NO$. The last assertion of the theorem
follows from Corollary \ref{11a}. \bs
 \begin{Cor}
For every $f\in\NO$ with the divisor $\Dab$ and $z_0\neq a_n,\,
z_0\neq b_n$ for all $n$  we have the representation
 \begin{equation}\label{rep1}
 f(z)=f(z_0)e^{\a (z-z_0)}\lim_{r\to\infty}{\prod_{n:|a_n|<r}(z-a_n)/
(z_0-a_n)e^{(z-z_0)/(a_n-z_0)}\over \prod_{n:|b_n|<r}(z-b_n)/
(z_0-b_n)e^{(z-z_0)/(b_n-z_0)}}.
 \end{equation}
 \end{Cor}
\bpr Note that for $|z_0|<r/4$ and $|z|<r/4$ we have
 \begin{equation}\label{AC}
\left|\log{\prod_{n:|a_n|<r}(1-z/a_n)e^{z/a_n}
\prod_{n:|b_n+z_0|<r}(1-z/b_n)e^{z/b_n}
\over\prod_{n:|a_n+z_0|<r}(1-z/a_n)e^{z/a_n}
\prod_{n:|b_n|<r}(1-z/b_n)e^{z/b_n}}\right|
 \end{equation}
 $$
\le\sum_{r-|z_0|\le |a_n|<r+|z_0|}|\log(1-z/a_n)+z/a_n|+
\sum_{r-|z_0|\le|b_n|<r+|z_0|}|\log(1-z/b_n)+z/b_n|
 $$
 $$
\le{C|z|^2\over(r-|z_0|)^2} [\card\{n:\,r-|z_0|\le|a_n|\le
r+|z_0|\}+\card\{n:\,r-|z_0|\le|b_n|\le r+|z_0|\}].
 $$
It follows from (\ref{d1}) that
$|\mu_D|(\{w:\,r-|z_0|<|w|<r+|z_0|\})=O(r)$ as $r\to\infty$. Hence
the right-hand side of (\ref{AC}) tends to $0$ as $r\to\infty$
uniformly on compacta in $\C$. Applying (\ref{rep}) with
$D=((a_n-z_0),(b_n-z_0))$ to $f(z+z_0)$, we obtain
(\ref{rep1}).\bs
\medskip

Using Theorem \ref{Y1} instead of \ref{Y}, we obtain the following
result:
\begin{Th}
Suppose a divisor $\Dab,\, a_n\neq 0,\,b_n\neq 0$ for all $n$,
satisfies conditions a), b), d), e), and f). Then $D$ is the
divisor of the normal function (\ref{rep}).
\end{Th}

\bigskip
\begin{center}
{\bf \S 3. Special properties of almost elliptic functions}
\end{center}
\bigskip

In what follows we need some  properties of almost periodic
mappings and divisors in the complex plain $\C$.
 \begin{Def}\label{rd}
A set $E\subset\C$ is called relatively dense if there exists
$L<\infty$ such that every disc of radius $L$ has a nonempty
intersection with $E$.
 \end{Def}
 \begin{Def}\label{ta}
Let $g$ be a continuous mapping  from $\C$ to a metric space
$(Y,\,d)$. A number $\tau\in\C$ is called an $\e$-almost period of
$g$ if
\begin{equation}\label{ap2}
d(g(z-\tau),g(z))<\e \quad \hbox{for all}\quad z\in\C.
\end{equation}
  The mapping  $g$ is called almost periodic if  for each
$\e>0$ the set of $\e$-almost periods of $g$ is relatively dense
in $\C$.
\end{Def}

The following results are well known for almost periodic functions
in the real axis (see, for example, \cite{B}, \cite{C}). One can
easily carry over their proofs to our case.

\begin{Pro}\label{0}
a) An almost periodic mapping is bounded and uniformly continuous,

b) if a sequence of almost periodic mappings converges uniformly
in $\C$, then its limit is also an almost periodic mapping.
\end{Pro}

\begin{Pro}\label{00}
 Suppose $f:\C\to\C$ is almost periodic function such that $\inf_\C|f(z)|>0$;
 then we have
 $$
f(z)=e^{g(z)+i(\b'x+\b''y)}, \quad \b',\b''\in\R,\quad z=x+iy,
 $$
where $g$ is an almost periodic function in $\C$.
\end{Pro}

 Furthermore, the following proposition is valid.
\begin{Pro}\label{1}
Suppose  $g:\C\to Y$ is a continuous mapping; then the following
conditions are equivalent:

a) $g$ is almost periodic,

b) for each $\e>0$ there exists $L<\infty$ such that every
interval $(a,\,a+L)$ of the real axis and every interval
$(ib,\,ib+iL)$ of the imaginary axis contains a point $\tau$
satisfying (\ref{ap2}),

c) for each sequence $(h_n)\subset\C$ there exists a subsequence
$(h_{n'})$ such that $d(g(z+h_{n'}),\,g(z+h_{m'}))\to 0$ as
$n',m'\to\infty$ uniformly on $\C$,

In addition,  if $(Y,\,d)$ is the plane $\C$ with the Euclidean
metric, then a) -- c) are equivalent to the  condition

d) there is a sequence of finite exponential sums
\begin{equation}\label{sum}
  S_k(z)=\sum_j c_{j,k} e^{i(\l_{j,k}x+\l'_{j,k}y)},
\quad \l_{j,k},\l'_{j,k}\in\R,\quad z=x+iy,
\end{equation}
 such that $S_k(z)-g(z)\to 0$
as $k\to\infty$ uniformly on $\C$.
\end{Pro}

\bpr The equivalence of a), c), and d) is well known for almost
periodic functions in the real axis $\R$ (see, for example,
\cite{B}, \cite{C}). In the same way, one can easily prove a
similar result in our case. Furthermore, the sum of two
$\e$-almost periods is an $2\e$-almost period, hence b) implies
that there is an $2\e$-almost period in every disc of radius
$\sqrt{2}L$. Therefore, b) implies a). On the other hand, let
$G:\,x\to g(x+iy)$ be the mapping from $\R$ to the space $\tilde
Y$ of continuous bounded functions $r(y),\,y\in\R,$ with the
distance $\tilde d(r_1,\,r_2)=\sup_{y\in\R} d(r_1(y),\,r_2(y))$.
Suppose that a mapping $g$ satisfies condition c); then for each
sequence $(h_n)\subset\R$ there is a subsequence $(h_{n'})$ such
that $\tilde d(G(x+h_{n'}),\,G(x+h_{m'}))\to 0$ as
$n',m'\to\infty$ uniformly in $\R$. Consequently, $G$ is an almost
periodic mapping. Hence for each $\e>0$ there exists $L'<\infty$
such that every interval $(a,\,a+L')\subset\R$ contains a point
$\tau$ with the property
 $$
 \sup_{z\in\C}d(g(z+\tau),g(z))=\sup_{x\in\R}\tilde d(G(x+\tau),\,G(x))<\e .
 $$
For the same reason, there exists $L''<\infty$ such that every
interval $(ib,\,ib+iL'')$ of the imaginary axis contains $\tau$
with property (\ref{ap2}). Hence b) is valid for
$L=\max\{L',L''\}$. \bs

\medskip
The class $\AE$ is just the set of all nonconstant meromorphic
almost periodic mappings from $\C$ to the  Riemann sphere. Hence
we have just proved Theorem SB.

\begin{Def}\label{D}
A number $\tau\in\C$ is an $\e$-almost period of a divisor $\Dab$
if there exist  bijections $\sigma:\,\N\to\N$ and
$\sigma':\,\N\to\N$ such that
 \begin{equation}\label{bi}
|a_n+\tau-a_{\sigma(n)}|<\e, \quad |b_n+\tau-b_{\sigma'(n)}|<\e
\quad \forall\quad n\in\N.
 \end{equation}
 The divisor is almost periodic if for each $\e>0$
there exists $L=L(\e)<\infty$ such that every disc $B(z,L)$
contains $\e$-almost period $\tau$. If, in addition,  (\ref{bi})
 holds with $\s(n)\equiv\s'(n)$, we say that the
divisor $D$ is almost periodic with a regular indexing.
\end{Def}
 For the case of positive divisors, Definition \ref{D} is very
close to the definition of an almost periodic zero set in a strip
(see \cite{L}, Appendix VI, \cite{T}, \cite{R}, and \cite{FRR}).

\medskip
{\bf Remark}. Note that all previous definitions and statements
are stable under any renumeration of $a_n$ and $b_n$. The same is
true for the property of a divisor to be almost periodic, because
such change means the replacing $\Dab$ by $\t
D=((a_{s(n)}),\,(b_{s'(n)}))$, where $s,\,s'$ are in general
different bijections $\N\to\N$. However, it is not hard to see
that in the general case the property of a divisor to have a
regular indexing may violate, although it survives in the case
$s=s'$).

\medskip
\begin{Pro}\label{3}
Suppose a divisor $\Dab$ is almost periodic; then (\ref{d}) holds.
If, in addition, $\Dab$ has a regular indexing, then
\begin{equation}\label{ind}
d_0= \sup_n|a_n-b_n|<\infty.
 \end{equation}
\end{Pro}
For almost periodic positive divisors in a strip this assertion
was obtained in \cite{FRR}.
\medskip

 \bpr  There is $L<\infty$ such that any disc $B(c,L)$ contains
 an 1-almost period $\tau$ of the divisor $D$. If (\ref{bi})
 holds for bijections $\s,\s'$, then  for any
$a_n\in\overline{B(c,1)}$ we get $a_{\s(n)}\in B(0,2+L)$. Hence,
 $$
 \card\{n:\,|a_n-c|\le 1\}\le \card\{n:\,a_n\in B(0,L+2)\}.
 $$
In the same way, we bound $\card\{n:\,|b_n-c|\le 1\}$.

Now suppose that the divisor $D$ has a regular indexing. Set
 $$
 M=\max_{n:|a_n|<L+1}|a_n-b_n|.
 $$
 For a term $a_{n'}$ take
an $1$-almost period $\tau\in B(-a_{n'},\,L)$. Using (\ref{bi})
with $\e=1$, we get  $|a_{\s(n')}|<L+1$. Hence,
$|b_{n'}-a_{n'}|<|b_{\s(n')}-a_{\s(n')}|+2\le M+2$. \bs

\begin{Pro}\label{5}
Let $D$ be an almost periodic divisor with a regular indexing;
then for every sequence $(h_k)\subset\C$ there is a subsequence
$(h'_k)\subset(h_k)$, a divisor $\wt D=((\t a_n)\,(\t b_n))$ with
a regular indexing, and bijections $\t\s(k,\cdot):\,\N\to\N$ such
that
\begin{equation}\label{bi1}
\sup_n|a_{\t\s(k,n)}+h'_k-\t a_n|\to 0,
\quad\sup_n|b_{\t\s(k,n)}+h'_k-\t b_n|\to 0\quad\hbox{as}\quad
k\to\infty.
 \end{equation}
\end{Pro}
 \bpr Take $\e>0$. By $E(\e/4)$ denote the union of discs
$B(\tau,\,\e/4)$ over all $\e/4$-almost periods of $D$. It follows
from Definition \ref{D} that each set $h_k+E(\e/4)$ intersects
with the disc $B(0,\,L(\e/4))$, moreover, the $m_2$-measure of the
intersection is at least $\pi(\e/4)^2/4$. We get
$$
 m_2(\bigcap_{m=1}^\infty\bigcup_{k=m}^\infty[B(0,L(\e/4))\cap
 (h_k+E(\e/4))])>0.
$$
Hence there is a point $z'$, which
 belongs to every set $h'_k+E(\e/4)$ for a subsequence
 $(h_{k'})\subset(h_k)$.  Therefore for each two
terms $h'_k,\,h'_s$ there exist $\e/4$-almost periods
$\tau,\,\tau'$ such that $|(h'_k-h'_s)-(\tau-\tau')|<\e/2$. Since
$\tau-\tau'$ is an $\e/2$-almost period of $D$, we obtain
 $$
\sup_n|a_n+h'_k-a_{\s(k,s,n)}-h'_s|<\e,\quad
\sup_n|b_n+h'_k-b_{\s(k,s,n)}-h'_s|<\e
 $$
 for some bijections $\s(k,s,\cdot):\,\N\to\N$. Using the diagonal
process and passing on to a subsequence if necessary, we get
 $$
\sup_n|a_n+h'_k-a_{\s(k,s,n)}-h'_s|<2^{-k},\quad
\sup_n|b_n+h'_k-b_{\s(k,s,n)}-h'_s|<2^{-k}\quad\forall\,k\in\N,\,s>k.
 $$
By definition, put
 $$
 \t\s(1,\cdot)=\s(1,2,\cdot),\,
\t\s(2,\cdot)=\s(2,3,\cdot)\circ\s(1,2,\cdot),\,
\t\s(3,\cdot)=\s(3,4,\cdot)\circ\s(2,3,\cdot)\circ\s(1,2,\cdot),\dots
 $$
Since
 $$
 |a_{\t\s(k,n)}+h'_k-a_{\t\s(k+1,n)}-h'_{k+1}|<2^{-k},\quad
 |b_{\t\s(k,n)}+h'_k-b_{\t\s(k+1,n)}-h'_{k+1}|<2^{-k}\quad\forall
 n,k\in\N,
 $$
 we see that there exist limits
 $$\t a_n=\lim_{k\to\infty}(a_{\t\s(k,n)}+h'_k),\qquad
 \t b_n=\lim_{k\to\infty}(b_{\t\s(k,n)}+h'_k),\qquad\forall n\in\N.
 $$
It can be easily checked that the divisor $\wt D=((\t a_n)\,(\t
b_n))$ is almost periodic, has a regular indexing, and satisfies
(\ref{bi1}).\bs

\begin{Th}\label{13}
The divisor $D$ of $f\in\AE$ is almost periodic. Moreover, for any
$\e>0$ there exists a relatively dense set of common $\e$-almost
periods of $f$ and $D$.
 \end{Th}
\bpr Take $\e>0$ and put $\d=\min\{\d_0,\,\e\}/C_0$, where $\d_0$
is from (\ref{sep}), and $C_0$ from (\ref{d}). Let $A$ be a
connected component of $A(\d)$, and $\wt A$ the union of $A$ and
all bounded connected components of $\C\setminus A$. Then
$\diam\wt A\le\e$ and $\wt A$ does not contain zeros and poles of
$f$ simultaneously. Now, by Proposition \ref{Y}, $1/C\le|f(z)|\le
C$ for $z\not\in A(\d)$. Clearly, there is $\eta=\eta(C)$ such
that for any $\eta$-almost period $\tau$ of $f$ we obtain
 $$
 1/(2C)\le|f(z+\tau)|\le 2C, \qquad
 |f(z+\tau)/f(z)-1|<1/2\qquad\forall z\not\in A(\d).
 $$
Hence the increment of $\arg f(z+\tau)$ along $\partial\wt A$
coincides with the increment of $\arg f(z)$, and the functions
$f(z+\tau)$ and $f(z)$ have the same numbers of zeros (or poles)
in $\wt A$. Consequently, there exist bijections $\s$ and $\s'$ of
$\N$ to $\N$ such that (\ref{bi}) holds and $\tau$ is a common
$\max\{\e,\eta\}$-almost period of $f$ and $D$. \bs

\begin{Th}\label{15}
 For each $\AE$-function with the divisor $\Dab$ there is an
indexing of zeros and poles with property (\ref{ind}).
 \end{Th}

\bpr The proof is based on a partition of the complex plane into
quadrilaterals subordinating to $\f\in\AE$. The idea of partition
belongs to F.~Sunyer-i-Balaguer \cite{S}, but his construction
contains a small inaccuracy (he does not consider the case when
the projections of the zero set and the pole set to the real and
imaginary axes are dense). We shall give a complete proof here.

Let $f$ be an $\AE$-function with the divisor $\Dab,\, a_n\neq
0,\,b_n\neq 0$ for all $n$. Set $r=\min\{|a_n|/4,|b_n|/4,\,
n=1,2,\dots\}$, and take $\d<\min\{r,\e_0,1\}/(2C_0)$ such that
 $$
 \d\neq|\Re a_n|,\,\d\neq|\Re b_n|,\,\d\neq|\Im a_n|,\,\d\neq|\Im
b_n|\quad\forall\,n,
 $$
 $$
 2\d\neq|a_n-a_k|,\,2\d\neq|b_n-b_k|,\,2\d\neq|a_n-b_k|\quad\forall\,n,k.
 $$
 Let $\cup_{k=1}^\infty A_k$ be the
decomposition of the set $A(\d)=\cup_n B(a_n,\d) \cup \cup_n
B(b_n,\d)$ into connected components. Note that
$\overline{A_k}\cap\overline{A_{k'}}=\emptyset$ for all $k\neq k'$
and any disc $B(a_n,\d)$ or $B(b_n,\d)$ is not tangent to the real
or imaginary axis.  Also, by Proposition \ref{i}, $\diam A_k<\e_0$
for all $k$, therefore each  $A_k$  does not contain zeros and
poles simultaneously.

Furthermore, let $A_{k_1}$ be the component with the minimal index
that intersects with the real axis $l$ and $(\a_1,\b_1)$ be the
minimal interval of $l$ containing this intersection.  Replace in
$l$ the interval $(\a_1,\b_1)$ by a Jordan curve
$L_1\subset\partial A_{k_1}$ with the same endpoints. Note that
the length of $L_1$ does not exceed $2\pi\d~\card\{n:\,a_n\hbox{
or }b_n\in A_{k_1}\}$. Take the component $A_{k_2}$ with the
minimal index that intersects with $l\setminus(\a_1,\b_1)$ and
repeat the above construction. Continuing like this, we obtain a
Jordan curve $l_x\subset\{z:\,|y|<r\},\,0\in l_x,$ such that
$l_x\cap A(\d)=\emptyset$.  Note that if $A_k$ intersects with a
segment $[\a,\,\b]$ of real axis, then $A_k$ is contained in the
rectangle $[\a-1,\,\b+1]\times[-1,\,1]$. In view of (\ref{d1}),
the number of terms $a_n,\,b_n$ in this rectangle is at most
$12\max\{\b-\a,1\}C_0$, therefore the length of $l_x$ inside the
rectangle $\{z:\,\a<x<\b,\,|y|<r\}$ is at most
$C'\max\{\b-\a,1\}$, where $C'$ depends on $\d$ and $\C_0$.

Similarly, we obtain a Jordan curve
$l_y\subset\{z:\,|x|<r\},\,0\in l_y,$ such that $l_y\cap
A(\d)=\emptyset$  and the length of $l_y$ inside any rectangle
$\{z:\,|x|<r,\,\a<y<\b\}$ is at most $C'\max\{\b-\a,1\}$.

In view of (\ref{pz}), there is $C_\d<\infty$ such that
 $$
 1/C_\d<|f(z)|<C_\d\qquad\forall\,z\not\in A(\d).
 $$
 Take $\e=\e(\d)$ such that
 \begin{equation}\label{arg}
 1/(2C_\d)\le|f(z+\tau)|\le 2C_\d, \qquad
 |f(z+\tau)/f(z)-1|<1/2,
  \end{equation}
whenever $\rho_S(f(z+\tau),\,f(z))<2\e$ and $z\not\in A(\d)$.
Using Definition \ref{S}, take $L>4r$ and pick up a sequence of
$\e$-almost periods $(\tau_p)_{p\in\Z}$ in the real axis and
another one $(i\tau'_q)_{q\in\Z}$ in the imaginary axis such that
 $$
L<\tau_{p+1}-\tau_p<3L,\quad L<\tau'_{q+1}-\tau'_q<3L.
 $$
Since $B(0,\,3r)\cap A(\d)=\es$, we see that (\ref{arg}) implies
$B(\tau,\,2r)\cap A(\d)=\es$. Taking into account embedding
$l_x\subset\{z:\,|y|<r\}$, $l_y\subset\{z:\,|x|<r\}$, we obtain
 $$
 B(\tau_p,2r)\cap l_x=B(\tau_p,2r)\cap\{z:\,\Im
 z=0\}\quad\forall\,p\in\Z,
 $$
 $$
 B(i\tau'_q,2r)\cap l_y=B(i\tau'_q,2r)\cap\{z:\,\Re
 z=0\}\quad\forall\,q\in\Z.
 $$
Hence,  the bounded domain $R_{p,q}$ formed by the lines
$l_x+i\tau'_q, \,l_y+\tau_p, \,l_x+i\tau'_{q+1},
\,l_y+\tau_{p+1}$,  is a quadrilateral. Using (\ref{arg}), we
obtain that the difference between the increments of continuous
branches of $\arg f(z)$ along opposite sides of each quadrilateral
$R_{p,q}$ is less than $\pi$. Now the Argument Principle yields
the equality
 $$
 \card\{n:\,a_n\in R_{p,q}\}=\card\{n:\,b_n\in R_{p,q}\}.
 $$
for all $p,\,q\in\Z$. Since $\diam R_{p,q}<15L$ for all $p,q$, we
obtain (\ref{ind}). \bs

\medskip
By definition, put
 $$
G_n(z)= \log\left({b_n(a_n-z)\over
a_n(b_n-z)}\right),\qquad n\in\N,\quad a_n,\,b_n\neq 0.
 $$
The function $G_n(z)$ is well defined in the complex plane with
discontinuity in $[a_n,\,b_n]$ under the condition $G_n(0)=0$.
Also, put $\Im G_n(z)=\pi+\Im G_n(\infty)$ for $z\in [a_n,\,b_n]$.
Suppose $f$ is an $\AE$-function with a divisor
$\Dab,\,a_n\neq0,\,b_n\neq0$.  By Theorem \ref{Rep}, $f$ has the
form (\ref{rep}). Using property (\ref{ind}), we obtain that the
function
 $$
  \arg f(z)=\Im\a z+\sum_n\left(\Im G_n(z)+\Im (z/a_n-z/b_n)\right)
 $$
is well defined in the complex plane, and so is $\log f(z)$
outside zeros and poles of $f$.

\begin{Pro}
Suppose  an $\AE$-function $f$ has the divisor $\Dab$ with a
regular indexing and property (\ref{ind}). Then there exists a
constant $\b\in\C$ such that
\begin{equation}\label{AA}
\arg f(z)=\Im(\b z) +O(1).
\end{equation}
\end{Pro}
\bpr Let $\f(z)$ be a nonnegative smooth function in $\C$ such
that $\f(z)=1$ if $\dist\{z,\,[0,1]\}<1/2$ and $\f(z)=0$ if
$\dist\{z,\,[0,1]\}>1$. By definition, put
\begin{equation}\label{H}
H(z)=\exp\left\{\sum_n \f\left({z-a_n\over
 b_n-a_n}\right)G_n(z)\right\}.
\end{equation}
Also, put $\d=\min\{1/2,\e_0\}/(5C_0)$. Take $\e>0$, and let
$\tau$ be a common $\e$-almost period of $f$ and $D$.
 It follows from
(\ref{d}) and (\ref{ind}) that every $z\in\C$ belongs to supports
of at most $C_0(2d_0+6)^2$ terms of the sum in (\ref{H}). Besides,
the divisor $D$ is  almost periodic with a regular indexing.
Therefore, $1/C<|H(z)|<C$ and $|H(z+\tau)-H(z)|<C\e$ for
$z\in\C\setminus[A(\d)\cup (A(\d)-\tau)]$. In view of (\ref{pz}),
we have $1/C_\d<|f(z)|<C_\d$ and $|f(z+\tau)-f(z)|<C\e$ as well.
Hence,
 \begin{equation}\label{H1}
1/C<\left|{f(z)\over H(z)}\right|<C,\quad \left|{f(z+\tau)\over
H(z+\tau)}-{f(z)\over H(z)}\right|<C\e, \quad \forall\,z\not\in
A(\d)\cup (A(\d)-\tau).
 \end{equation}
 Let $A$ be a connected component of
$A(\d)$. Since $\diam A<\min\{1/2,\e_0\}$, we see that $A$
contains either zeros, or poles of $f$. For example, suppose that
$A$ contains zeros of $f$. By definition, put
 $$
 H_A(z)=\exp\left\{\sum_{n:a_n\in A}\f\left({z-a_n\over
 b_n-a_n}\right)G_n(z)\right\}.
 $$
 We see that $H_A(z)$ is holomorphic and has the same zeros as $f(z)$
 in $A$. On the other hand, $1/C<|H(z)/H_A(z)|<C$ for
 $z\in\overline{A}$. Applying the Maximum Modulus Principle to the function $f(z)/H_A(z)$,
 we get the bounds (\ref{H1}) for $z\in A$,
 therefore, for all $z\in\C$.

Thus, the function $F(z)/H(z)$ is almost periodic in $\C$.
Applying Proposition \ref{00}, we get
 $$
 f(z)/H(z)=\exp\{g(z)+i\Im(\b z)\}\quad\hbox{for some}\quad
 \b\in\C.
 $$
Since the functions $g(z)$ and $\Im[\sum_n \f({z-a_n\over
 a_n-b_n})G_n(z)]$ are uniformly bounded in $\C$, we obtain (\ref{AA}).\bs

\medskip
For an arbitrary divisor $D$ with properties (\ref{d}) and
(\ref{ind}) denote by $\nu_D$ the discrete measure with the
complex masses $a_n-b_n$ at the points $c_n=(a_n+b_n)/2$. Clearly,
\begin{equation}\label{o}
  |\nu_D|(\overline{B(c,1)})<C\quad\forall\, c\in\C.
\end{equation}

\begin{Pro}\label{7}
Suppose $f$ is an $\AE$-function, $\{R_{p,q}\}$ is the above
partition of the plane into quadrilaterals, and the divisor $\Dab$
of $f$ has a regular indexing such that $\{n:\,a_n\in
R_{p,q}\}=\{n:\,b_n\in R_{p,q}\}$. Then for any simply connected
bounded domain $E\subset\C$  we have
 \begin{equation}\label{reg}
   |\nu_D(E)-\bar\b m_2(E)/2\pi|\le C(m_1(\partial E)+1).
 \end{equation}
\end{Pro}

\bpr Set
 $$
 I_1=\{(p,q):\,\overline{R}_{p,q}\cap E\ne\es\},\quad
I_2=\{(p,q):\,\overline{R}_{p,q}\subset E\},\quad R=\cup_{(p,q)\in
I_2}\overline{R}_{p,q}.
 $$
It follows from the definition of $R_{p,q}$ that for every $p,q$
there is $z_{p,q}$ such that $B(z_{p,q},\,1)\subset R_{p,q}\subset
B(z_{p,q},\,4L)$ and each point $z\in\C$ is contained in at most
$C(L)$ discs $B(z_{p,q},\,5L)$. Since $\partial E$ is connected,
we obtain that either $\partial E$ is contained in the unique disc
$B(z_{p,q},\,5L)$, or $m_1(\partial E\cap B(z_{p,q},\,5L))\ge L$
for all $(p,\,q)\in I_1\setminus I_2$. In both cases we obtain
 $$
 \card(I_1\setminus I_2)\le C(m_1(\partial E)+1),\quad
 m_1(\partial R)\le C(m_1(\partial E)+1),
 $$
\begin{equation}\label{a1}
  (m_2+|\nu_D|)(E\setminus R)\le(m_2+|\nu_D|)(\cup_{(p,\,q)\in I_1\setminus I_2}\overline{R}_{p,q})
  \le C(m_1(\partial E)+1).
\end{equation}
Furthermore, using (\ref{ind}), we get
 $$
 \card(\{n:\,c_n\in R,\,a_n,\,b_n\not\in
R\}\cup \{n:\,c_n\not\in R,\,a_n,\,b_n\in R\})\le Cm_1(\partial
R).
 $$
Therefore,
\begin{equation}\label{a}
\left|\nu_D(R)-\sum_{a_n,b_n\in R}(a_n-b_n)\right|\le
Cm_1(\partial R).
\end{equation}
Define a continuous branch $\widetilde{\arg} f(z)$ of the argument
 $f$ in the set $\cup_{p,q}\partial R_{p,q}$ by the condition
$\widetilde{\arg} f(0)=\arg f(0)$. The increment of the argument
of $f$ along $\partial R_{p,q}$ equals $0$ for all $p,\,q$,
therefore the branch is well defined. The both ends of each
segment $[a_n,\,b_n]$ belong to the same quadrilateral, hence the
sum of jumps of $\arg f(z)$ along each side of $R_{p,q}$ is zero.
Therefore,  $\widetilde{\arg} f(\tau_p+i\tau'_q)=\arg
f(\tau_p+i\tau'_q)$ for all $p,q$.
 Note that $|(\arg f)'(z)|=|(\widetilde{\arg} f'(z)|\le|f'(z)/f(z)|<C$ in $\cup_{p,q}\partial
 R_{p,q}$. Using (\ref{AA}), we get
 $$
 |\widetilde{\arg} f(z)-\Im(\b z)|<C\quad \forall z\in\cup_{p,q}\partial
 R_{p,q}.
 $$
Integrating by parts, we obtain
 $$
\sum_{a_n,b_n\in R}(a_n-b_n)={1\over 2\pi i}\int_{\partial
R}z{f'(z)\over f(z)}~dz=-{1\over 2\pi} \int_{\partial
R}\widetilde\arg f(z)~dz+O(m_1(\partial R))
 $$
 $$
=-{1\over 2\pi} \int_{\partial R}\Im(\b z)~dz +O(m_1(\partial
R))=\bar\b m_2(R)/2\pi+O(m_1(\partial R)).
 $$
Now (\ref{a1}) and (\ref{a}) yield the assertion of the
proposition.\bs

\bigskip
\begin{center}
{\bf \S 4. Almost elliptic functions with divisors having regular
indexing.}
\end{center}
\bigskip

Let $\mu^\tau$ be the translation of a measure $\mu$ in $\C$,
i.e., $\mu^\tau(E)=\mu(E+\tau)$ for each Borel set $E\subset\C$.

Here we shall prove the following theorem:
\begin{Th}\label{main}
 An  almost periodic divisor $\Dab,\,
a_n\neq 0,\,b_n\neq 0$ for all $n$, with a regular indexing is the
divisor of $f\in\AE$ if and only if the following conditions are
fulfilled:
\medskip

a) the divisor $D$ obeys (\ref{sep}),
\medskip

b) the measure $\nu_D$ satisfies (\ref{reg}) for every convex
bounded subset of $\C$,
\medskip

c) the measure $\nu_D$ satisfies the bound \footnote{It follows
from (\ref{l5}) that we can check this condition only at points of
a relatively dense subset.}
\begin{equation}\label{h2}
\limsup_{r\to\infty}\left|\int_{1\le |w|\le
r}{d\nu_D^z(w)-d\nu_D(w)\over w}\right|<C\qquad\forall z\in\C.
 \end{equation}
In this case the function $f$ has the form
\begin{equation}\label{f1}
f(z)=e^{\b z/2}\lim_{r\to\infty}\prod_{n:|c_n|<r}e^{G_n(z)}
 \end{equation}
up to a constant factor.
\end{Th}

 We begin with some lemmas:

\begin{Lem}
Suppose an  almost periodic divisor $\Dab$ has a regular indexing;
then we have
\begin{equation}\label{b1}
\left|\int_{E^+(r)}(1/w-{\overline w}/
r^2)d\nu(w)\right|+\left|\int_{E^-(r)}(1/w-{\overline w}/
r^2)d\nu(w)\right|<C{1+|z|^2\over r},\quad r>1,
\end{equation}
where
 $$
E^+(r)=\{w:\,|w|\le r,\,|w+z|>r\},\qquad E^-(r)=\{w:\,|w+z|\le
r,\,|w|>r\}.
 $$

Moreover, if b) holds for $\nu_D$ and $\l=\nu_D-(\bar\b/2\pi)
m_2$, then we have
\begin{equation}\label{l2}
\left|\int_{|w|\le r}w\,d\l(w)\right|\le Cr^2,\quad
\left|\int_{|w|\le r}\overline{w}\,d\l(w)\right|\le Cr^2,\quad
r>1,
\end{equation}
\begin{equation}\label{l3}
\left|\int_{1\le |w|\le r}d\l(w)/w\right|\le C(1+\log r),\quad
r>1,
\end{equation}
\begin{equation}\label{l4}
\left|\int_{r\le |w|\le R}d\l(w)/w^2\right|\le Cr^{-1},\quad
R>r>1.
\end{equation}
\begin{equation}\label{l5}
\left|\int_{1\le |w|\le r}{d\l^z(w)-d\l(w)\over w}\right|\le
C(1+\log^+|z|),\quad r>|z|^2+1,
\end{equation}
\begin{equation}\label{l6}
\left|\int_{1<|w|\le r, 1<|w-z|}\left({1\over w-z}-{1\over
 w}\right)d\nu_D(w)-\int_{1\le
|w|\le r}{d\nu_D^z(w)-d\nu_D(w)\over w}+{\bar\b\bar z\over
2}\right|<C.
\end{equation}
 All constants do not depend on shifts of $D$.
\end{Lem}

\bpr Using Proposition \ref{3}, we obtain (\ref{d}) and
(\ref{ind}). Hence the measure $\nu$ is well defined and satisfies
(\ref{o}). Hence, we obtain the bound $|\nu|(E^+(r)\cup
E^-(r))<C(1+|z|)r$. Since $|1/w-\overline{w}/r^2|<C|z|(r+|z|)/r^3$
for $w\in(E^+(r)\cup E^-(r))$, we get (\ref{b1}).

Next, put $w=u+iv$, $\a_1(t)=\l\{w:\,|w|\le r,\,u<t\}$.  We have
 \begin{equation}\label{C}
\int_{|w|\le r}u\,d\l(w)=\int_{-r}^r
t\,d\a_1(t)=r\a_1(r)-\int_{-r}^r \a_1(t)\,dt.
 \end{equation}
Since $|\a_1(t)|\le Cr$, we see that the module of (\ref{C}) has
the bound $Cr^2$. Clearly, $\int_{|w|\le r}v\,d\l(w)$ has the same
bound $Cr^2$. Hence, we get (\ref{l2}).

Furthermore, integrating by parts, we get
  $$
\int\limits_{1\le |w|\le r}{d\l(w)\over w}=\int\limits_{1<|w|\le
r}{\overline{w}~d\l(w)\over |w|^2}={1\over
r^2}\int\limits_{1<|w|\le
r}\overline{w}~d\l(w)+2\int\limits_1^r\left(\int\limits_{1<|w|\le
t}\overline{w}~d\l(w)\right){dt\over t^3}.
 $$
 Now (\ref{l2}) implies (\ref{l3}).

To prove (\ref{l4}), consider the integral
  \begin{equation}\label{tem}
\int_{r<u^2+v^2\le R, u>0}u^2/(u^2+v^2)^2d\l(w)
 =\int_r^\infty(1/t^2)d\a_2(t)
 =2\int_r^\infty\a_2(t)/t^3\,dt,
 \end{equation}
 where $\a_2(t)=\l(B(t/2,\,t/2)\cap[B(0,\,R)\setminus B(0,\,r)])$.
 Using (\ref{reg}), we get
  $$
  |\a_2(t)|\le|\l(B(0,\,R)\cap B(t/2,\,t/2))|+|\l(B(0,\,r)\cap B(t/2,\,t/2))|\le Ct.
  $$
 Hence the modulus of integral (\ref{tem}) does not
 exceed $C/r$. Clearly,  the same bound is valid for the integral
\begin{equation}\label{tem1}
\int\limits_{r<|w|\le R}\Re(1/w^2)d\l(w)=\int\limits_{r<|w|\le
R}u^2/(u^2+v^2)^2 d\l(w)-\int\limits_{r<|w|\le R}v^2/(u^2+v^2)^2
d\l(w).
  \end{equation}
 The orthogonal transformation of coordinates
 $u'=(u+v)/\sqrt{2},\,v'=(u-v)/\sqrt{2}$ reduces the integral
 $$
\int_{r<|w|\le R}\Im(1/w^2)d\l(w)= \int_{r<|w|\le
R}(2uv)/(u^2+v^2)^2 d\l(w)
 $$
 to (\ref{tem1}), so (\ref{l4}) follows.

 To prove bound (\ref{l5}), put $r_1=(|z|+1)^2$.
 Decompose integral in (\ref{l5}) into the sum
 $$
\int_{1\le |w|\le r_1}{d\l^z(w)-d\l(w)\over w} +z\int_{r_1<|w|\le
r}{d\l(w)\over w^2} +z^2\int_{r_1<|w|\le r}{d\l(w)\over w^2(w-z)}+
\int_{E^-(r_1)}{d\l^z(w)\over w}
 $$
 $$
-\int_{E^+(r_1)}{d\l^z(w)\over w}+ \int_{E^+(r)}{d\l^z(w)\over w}
-\int_{E^-(r)}{d\l^z(w)\over w}
 =I_1+I_2+I_3+I_4-I_5+I_6-I_7.
 $$
 In view of (\ref{l3}) and (\ref{l4}), we have $|I_1|<C(1+\log r_1)$ and
$|I_2|<C|z|/r_1$. Using (\ref{m2}) and inequality $|z|<r_1/2\le
t/2$, we get
 $$
 |I_3|<2|z|^2\int_{r_1}^r{d|\l|(B(0,\,t))\over t^3}<
 C|z|^2/r_1.
 $$
Next,
 $$
I_6-I_7=\int_{E^+(r)}\left({1\over w}-{{\overline w}\over
r^2}\right)d\l^z(w)-\int_{E^-(r)}\left({1\over w}-{{\overline
w}\over r^2}\right)d\l^z(w)+\int_{|w|\le r}{{\overline w}\over
r^2}d\l^z(w)-
 $$
 $$
 + \int_{|w|\le r}{{\overline w}\over
r^2}d\l^z(w))+{{\overline z}\l^z(B(0,\,r))\over
r^2}=I_8-I_9+I_{10}-I_{11}+I_{12}.
 $$
Applying (\ref{b1}) with the measure $\l^z$, we get
$|I_8|+|I_9|<C$. Taking into account (\ref{l2}) and (\ref{reg}),
we get $|I_{10}|<C$, $|I_{11}|<C$, and $|I_{12}|<C|z|/r<C$.  In
the same way, $|I_4-I_5|<C$. So (\ref{l5}) is proved.

 To prove (\ref{l6}), note that
 $$
\left|\int_{\overline{B(z,1)}\setminus B(0,1)}{d\nu_D(w)\over
w}\right|+\left|\int_{\overline{B(0,1)}\setminus
B(z,1)}{d\nu_D(w)\over w-z}\right|<C/(1+|z|).
 $$
Next,
 $$
  \int_{1<|w|\le r,1<|w-z|}{d\nu_D(w)\over w-z}=
  \int_{1<|w|\le r,1<|w+z|}{d\nu^z_D(w)\over w}+
  {1\over r^2}\int_{|w+z|\le r}\overline{w}d\nu^z_D(w)
 $$
 $$
- {1\over r^2}\int_{|w|\le r}\overline{w}d\nu^z_D(w) +
\int_{E^+(r)}\left({1\over w}-{\overline{w}\over
  r^2}\right)d\nu^z_D(w)-\int_{E^-(r)}\left({1\over
w}-{\overline{w}\over r^2}\right)d\nu^z_D(w).
 $$
In view of (\ref{b1}), the difference of last two integrals is
uniformly bounded. Furthermore,
 $$
 \int_{|w+z|\le r}\overline{w}d\nu^z_D(w)-\int_{|w|\le r}\overline{w}d\nu^z_D(w)=
\int_{|w|\le r}\overline{w}d\l(w)-\int_{|w|\le
r}\overline{w}d\l^z(w)
 $$
 $$
-\overline{z}\left(\int_{|w|\le r}d\l(w)+ {\bar\b\over
2\pi}\int_{|w|\le r}dm_2(w)\right).
 $$
Applying (\ref{l2}) to the measures $\l$ and $\l^z$, and
 (\ref{reg}) to $\l(B(0,\,r))$, we obtain (\ref{l6}). The proof
 is complete. \bs

\begin{Lem}\label{e}
 Suppose an almost periodic divisor $\Dab,\,a_n\neq0,\,b_n\neq0,$
 with a regular indexing satisfies b).  Then the limit
 \begin{equation}\label{L}
G(z)=\lim_{r\to\infty}\sum_{n:|c_n|<r}G_n(z)
 \end{equation}
 exists uniformly on compacta in $\C$. Moreover, for every $\d>0$
 there exists a constant $C_\d$ such that the inequality
\begin{equation}\label{h1}
\limsup_{r\to\infty}\left|G(z)+ {\bar\b\bar z\over 2}-\int_{1\le
|w|\le r}{d\nu_D^z(w)-d\nu_D(w)\over w}\right|<C_\d
\end{equation}
is valid  for all $z\in \C\setminus A(\d)$.
 \end{Lem}

\bpr If $|c_n-z|>3d_0$ and $|c_n|>3d_0$,  then we get
 $$
G_n(z)=\log{a_n-z\over b_n-z}-\log{a_n\over b_n},
 $$
 $$
\left| \log\left({a_n-z\over b_n-z}\right)-{a_n-b_n\over
b_n-z}+{(a_n-b_n)^2\over2(b_n-z)^2}\right|<C\left(\left|{a_n-b_n\over
b_n-z}\right|^3\right),
 $$
 and
 $$
{a_n-b_n\over b_n-z}-{(a_n-b_n)^2\over2(b_n-z)^2}=
{\nu_D(\{c_n\})\over c_n-z}-{(a_n-b_n)^3\over4(b_n-z)^2(c_n-z)}.
 $$
Therefore,
 \begin{equation}\label{log1}
|G_n(z)-\nu_D(\{c_n\})/(c_n-z)+\nu_D(\{c_n\})/
c_n|<C(|a_n-z|^{-3}+|a_n|^{-3}).
 \end{equation}

Take $R>r>2|z|$. Proposition \ref{3} implies (\ref{d}). Therefore
the sum $\sum_{|a_n-z|>1}|a_n-z|^{-3}$ converges uniformly in
$\C$. Hence the sum $\sum_{n:r<|c_n|\le R}G_n(z)$  equals the
integral
 \begin{equation}\label{h}
\int_{r<|w|\le R}\left({1\over w-z}-
 {1\over w}\right)d\nu_D(w)=z\int_{r<|w|\le R}{d\nu_D(w)\over w^2}
+z^2\int_{r<|w|\le R}{d\nu_D(w)\over w^2(w-z)}
 \end{equation}
up to the term tending to zero as $r\to\infty$.  Taking into
account (\ref{l4}) and equality $\int_{r<|w|\le R}dm_2/w^2=0$, we
see that the first integral in the right-hand side of (\ref{h})
tends to $0$ as $r\to\infty$. Using (\ref{m2}), we get
 $$ \left|\int_{r<|w|\le
R}{d\nu_D(w)\over w^2(w-z)}\right|<
\int_r^R{d\,|\nu_D|(B(0,t))\over(t-|z|) t^2}\to 0
 $$
as $r\to\infty$. Hence, the limit in (\ref{L}) exists.

Thus, using (\ref{log1}), we  get
 \begin{equation}\label{s}
\left|\sum_{n:3d_0<|c_n|<r,|c_n-z|>3d_0}G_n(z)-
 \int_{3d_0<|w|\le r, 3d_0<|w-z|}\left({1\over w-z}-{1\over
 w}\right)d\nu_D(w)\right|<C.
 \end{equation}
Taking into account (\ref{o}), replace $3d_0$ by $1$.  Combining
(\ref{l6}) and (\ref{s}), we get (\ref{h1}).\bs

\medskip
\noindent {\bf Proof of the Theorem. Necessity}.

It follows from Theorem \ref{8} and Propositions \ref{3}, \ref{7}
that the divisor $\Dab$ of $\AE$-function $f$ is almost periodic
and satisfies  conditions (\ref{sep}) and (\ref{d}). Moreover,
since $D$ has a regular indexation,  we get b) and (\ref{ind}).
Next, for every $\d>0$ bounds (\ref{pz}) and (\ref{AA}) imply the
estimate
 \begin{equation}\label{b3}
|\log|f(z)|+i\arg f(z)-i\Im(\b
z)|<C_\d\qquad\forall\,z\in\C\setminus A(\d).
 \end{equation}
By Lemma \ref{e}, the function $G(z)$ is well defined. Taking into
account the definition of $\arg f(z)$, we see that the function
 $$
 E(z)=\log f(z)-G(z)-\b z/2
 $$
 is holomorphic in $\C$. Using (\ref{h1}) and (\ref{b3}), we get
 $$
\limsup_{r\to\infty}\left|E(z)+\int_{1\le |w|\le
r}{d\nu_D^z(w)-d\nu_D(w)\over w}\right|<C'_\d
 $$
for a sufficiently small $\d$ and  $z\in\C\setminus A(\d)$. It
follows from (\ref{l5}) that $|E(z)|<C(1+\log^+|z|)$ in this set.
Therefore,  $E(z)\equiv\const$, and we obtain (\ref{f1}) and
(\ref{h2}).
\medskip

\noindent {\bf Sufficiency}. Let $\Dab$ be a divisor satisfying
conditions of the Theorem. Let us show that the function $f$
(\ref{f1}) belongs to the class $\AE$.

  First, it follows from (\ref{h1}) and (\ref{h2}) that
$f$ satisfies (\ref{pz}). The divisor $D$ is almost periodic with
a regular indexing, hence Proposition \ref{3} implies (\ref{d})
and (\ref{N1}). By Theorem \ref{Y}, we get $f\in\NO$.

Take $\d_1=\min\{\d_0/6C_0,1/3\}$. Put
$A(D,\d_1)=\cup_nB(a_n,\d_1)\cup\cup_nB(b_n,\d_1)$ and
 $$
\a_1(f)= \sup\{\log|f(z)|:\,z\in\C\setminus A(D,\d_1)\},
 $$
 $$
  \a_2(f)= \inf\{\log|f(z)|:\,z\in\C\setminus A(D,\d_1)\},
 $$
Note that all zeros of $f$ belong to the same connected component
of $\C\setminus\cup_n B(b_n,\,\d_1)$, and all poles of $f$ belong
to the same connected component of $\C\setminus\cup_n
B(a_n,\,\d_1)$. The function $g(z)=\Im G(z)- \Im\b z/2$ is
harmonic in $\C\setminus\cup_n[a_n,\,b_n]$ and has the jump $2\pi$
in each segment $[a_n,\,b_n]$. It follows from (\ref{h1}) and
(\ref{h2}) that $g(z)$ is uniformly bounded in $\C$. Put
 $$
\a_3(g)=\sup\{g(z):\,z\in\C\},
 $$
 $$
\a_4(g)=\inf\{g(z):\,z\in\C\}.
 $$
Let $(h_k)$ be an arbitrary sequence.   It follows from
Proposition \ref{5} that there are an almost periodic divisor with
a regular indexing $\wt{D}=((\t a_n),(\t b_n))$ and a subsequence
$(h'_k)$, which satisfy (\ref{bi1}). Passing on to a subsequence,
we can suppose that the functions $g(z-h'_k)$ converge uniformly
on compacta in $\C\setminus\cup_n[\t a_n,\,\t b_n]$ to a function
$\t g$. Clearly, $\t g$ is harmonic in $\C\setminus\cup_n[\t
a_n,\,\t b_n]$ and has the jump $2\pi$ in each segment $[\t
a_n,\,\t b_n]$. Moreover, we have
\begin{equation}\label{lim}
 \a_3(\t g)\le\a_3(g),\quad
\a_4(\t g)\ge\a_4(g).
\end{equation}
Furthermore, since $f\in\NO$, we can pass on to a subsequence to
obtain
 $$
\rho_S(f(z-h'_k),\,\t f(z))\to 0\quad \hbox{as}\quad k\to\infty
 $$
 uniformly on compacta in $\C$. Clearly, $\t f(z)$ has the divisor $\wt{D}$
and the sequence $(\log|f(z-h'_k)|)$ converges to $\log|\t f(z)|$
uniformly on compacta in $\C\setminus A(\wt{D},\d_1)$. Hence,
\begin{equation}\label{lim2}
 \a_1(\t f)\le\a_1(f),\quad
\a_2(\t f)\ge\a_2(f).
\end{equation}
Since  $f(z)=|f(z)|\exp\{ig(z)+i\Im\b z\}$, we get
  \begin{equation}\label{f}
\t f(z)=|\t f(z)|e^{i\t g(z)+i\Im\b
z}e^{i\t\th},\quad\hbox{where}\quad e^{i\t\th}=
\lim_{k\to\infty}e^{-i\Im(\b h'_k)}.
 \end{equation}

Let us check that inequalities (\ref{lim}) and (\ref{lim2}) are in
fact equalities.  Rewrite (\ref{bi1}) in the form
 $$
\sup_n|\t a_{\s^{-1}(k,n)}-h'_k-a_n|\to 0,\quad \sup_n|\t
b_{\s^{-1}(k,n)}-h'_k-b_n|\to 0 \quad\hbox{as}\quad k\to\infty.
 $$
For some subsequence $(h''_k)\subset (h'_k)$ the functions $\t
g(z+h''_k)$ converge uniformly on compacta in
$\C\setminus\cup_n[a_n,\,b_n]$ to a function $\widehat{g}$, which
is harmonic in $\C\setminus\cup_n[a_n,\,b_n]$ and has the jump
$2\pi$ in each segment $[a_n,\,b_n]$. As above, for some
$\NO$-function $\widehat{f}$ with the divisor $D$ we have
$$ \rho_S(\t f(z+h''_k),\,\widehat{f}(z))\to 0\quad \hbox{as}\quad
k\to\infty
 $$
 uniformly on compacta in $\C$. Therefore, the function $f(z)/\widehat
 {f}(z)$ is an entire function without zeros, which satisfies the
 inequality $1/C\le |f(z)/\widehat{f}(z)|\le C$ in the set $\C\setminus
 A(D, \d_1)$. Hence, $\widehat{f}(z)\equiv Kf(z),\,K\in\C$. Now the
 inequalities
 $$
  \a_1(\widehat{f})\le\a_1(\t f)\le\a_1(f),\quad
\a_2(\widehat{f})\ge\a_2(\t f)\ge\a_2(f)
 $$
imply that $|K|=1$, and equalities (\ref{lim2}) prevail.
Furthermore, using (\ref{f}), we get
 $$
\widehat{f}(z)=|\widehat{f}(z)|e^{i\widehat{g}(z)+i\Im\b
z}e^{i\t\th}e^{i\widehat{\th}},\quad\hbox{where}\quad
e^{i\widehat{\th}}= \lim_{k\to\infty}e^{i\Im(\b h''_k)}.
 $$
Since $e^{i\widehat{\th}}=e^{-i\t\th}$ and $\widehat{f}=e^{i\arg
K}f$, we get $\widehat{g}(z)=\arg K+2\pi l+g(z),\,l\in\Z$. Hence
the inequalities
 $$
 \a_3(\widehat{g})\le \a_3(\hat g)\le\a_3(g),\quad
\a_4(\hat g)\ge\a_4(\t g)\ge\a_4(g),
 $$
imply that $K=1$, and equalities (\ref{lim}) prevail.

  To prove $f\in\AE$  we shall show  that
$\sup_{z\in\C}\rho_S(f(z-h'_k),\t f(z))\to 0$ as $k\to\infty$.

Assume the contrary. Then there exists an $\e_0>0$ and a sequence
$(z_k)$ such that
\begin{equation}\label{cont}
\rho_S(f(z_k-h'_k),\t f(z_k))>\e_0.
\end{equation}
Using Proposition \ref{5}, take a subsequence $(t_{k'})$ of
$(h''_k)$ with the following properties:
 \begin{equation}\label{bi3}
 \sup_n|a_{\s^*(k',n)}+t_{k'}-z_{k'}-a^*_n|\to 0,\quad
 \sup_n|b_{\s^*(k',n)}+t_{k'}-z_{k'}-b^*_n|\to 0 \quad\hbox{as}\quad
k'\to\infty,
  \end{equation}
\begin{equation}\label{bi4}
 \sup_n|\t a_{\s^{**}(k',n)}-z_{k'}-a^{**}_n|\to 0,\quad
 \sup_n|\t b_{\s^{**}(k',n)}-z_{k'}-b^{**}_n|\to 0 \quad\hbox{as}\quad
k'\to\infty,
  \end{equation}
 $$
 g(z-t_{k'}+z_{k'})\to g^*(z)\quad\hbox{as}\quad k'\to\infty
 $$
 uniformly on compacta in $\C\setminus\cup_n[a_n^*,\,b_n^*]$,
 $$
 \t g(z+z_{k'})\to g^{**}(z)\quad\hbox{as}\quad k'\to\infty
 $$
 uniformly on compacta in
 $\C\setminus\cup_n[a_n^{**},\,b_n^{**}]$,
 \begin{equation}\label{f4}
\rho_S(f(z-t_{k'}+z_{k'}),\,f^*(z))\to 0,\quad\rho_S(\t
 f(z+z_{k'}),\,f^{**}(z))\to 0\quad\hbox{as}\quad k'\to\infty
 \end{equation}
 uniformly on compacta in $\C$.

  Here $D^*=((a_n^*),\,(b_n^*)),\,
D^{**}=((a_n^{**}),\,(b_n^{**}))$, $\s^*(k',\cdot)$ and
$\s^{**}(k',\cdot)$ are bijections $\N\to\N$, $f^*$ is an
$\NO$-function in $\C$ with the divisor $D^*$, $f^{**}$ is an
$\NO$-function in $\C$ with the divisor $D^{**}$, $g^*$ is a
harmonic function in $\C\setminus\cup_n[a_n^*,\,b_n^*]$, which has
the jump $2\pi$ in each segment $[a_n^*,\,b_n^*]$, and $g^{**}$ is
a harmonic function in $\C\setminus\cup_n[a_n^{**},\,b_n^{**}]$,
which has the jump $2\pi$ in each segment $[a_n^{**},\,b_n^{**}]$.

Combining (\ref{bi1}), (\ref{bi3}), and (\ref{bi4}), we get
 $$
\sup_n|a^*_{\widehat{\s}(k',n)}-a^{**}_n| \to 0, \quad \sup_n
|b^*_{\widehat{\s}(k',n)}-b^{**}_n|\to 0\quad\hbox{as}\quad
k'\to\infty,
 $$
where
$\widehat{\s}(k',\cdot)=(\s^*)^{-1}(k',\cdot)\circ\s(k',\cdot)\circ\s^{**}(k',\cdot)$.
The sequence $(a^*_n)$ has not limit points, hence
 for each $n$ the number $\widehat{\s}(k',n)$ is the same for all $k'>k_n$.
Since $\widehat{\s}(k',n)$ are bijections, we obtain that
$\widehat{\s}(n)=\lim_{k'\to\infty}\widehat{\s}(k',n)$ is a
bijection as well. Therefore,
 $$
 a^*_{\widehat{\s}(n)}=a^{**}_n,\quad
 b^*_{\widehat{\s}(n)}=b^{**}_n\quad\forall\,n.
 $$
 Hence, $D^{**}=D^*$ up to the same rearrangements of $(a^*_n)$
  and $(b^*_n)$.

Furthermore, (\ref{f1}) and (\ref{f}) yield
 $$
f^*(z)=|f^*(z)|e^{ig^*(z)+i\Im\b
z}e^{i\th^*},\quad\hbox{where}\quad e^{i\th^*}=
\lim_{k\to\infty}e^{-i\Im\b(t_{k'}-z_{k'})},
 $$
 $$
f^{**}(z)=|f^{**}(z)|e^{ig^{**}(z)+i\Im\b
z}e^{i\t\th}e^{i\th^{**}},\quad\hbox{where}\quad e^{i\th^{**}}=
\lim_{k\to\infty}e^{i\Im\b z_{k'}}.
 $$
 As before, the function
$f^*(z)/f^{**}(z)$ is bounded on the set $\C\setminus A(D^*,\d_1)$
and has not zeros and poles in $\C$. Hence, $f^*(z)=Kf^{**}(z)$.
The equalities $\a_j(f^*,D^*)=\a_j(f,D)=\a_j(\t
f,\wt{D})=\a_j(f^{**},D^{**}),\,j=1,2$ show that $|K|=1$. Since
$e^{i\th^*}=e^{i\t\th}e^{i\th^{**}}$, we get
$g^*(z)=g^{**}(z)+\arg K+2\pi l,\,l\in\Z$.  The equalities
$\a_j(g^*,D^*)=\a_j(g,D)=\a_j(\t
g,\wt{D})=\a_j(g^{**},D^{**}),\,j=3,4$ yield $g^*(z)\equiv
g^{**}(z)$ and $f^*(z)\equiv f^{**}(z)$. The last equality
contradicts to (\ref{f4}) and (\ref{cont}). Hence, $f\in\AE$. \bs

\bigskip
\begin{center}
{\bf \S 5. Examples of almost elliptic functions}
\end{center}
\bigskip

Suppose that $Q(z)$ is a finite exponential sum (\ref{sum}) (or a
uniform limit of a sequence of sums (\ref{sum})) with the
additional condition
\begin{equation}\label{sp}
\l_{j,k}n'+\l'_{j,k}n''\not\in 2\pi\Z\quad\hbox{for all}\quad
n=n'+in''\in\Z^2.
\end{equation}
Consider the almost periodic mapping $F(z)=(Q(z);\exp2\pi
ix;\exp2\pi iy)$ from $\C$ to $\C^3$ with the Euclidean metric. It
can  be easily checked that $F(z)$ satisfies condition b) of
Proposition \ref{1}. Therefore, $F(z)$ is almost periodic. Then,
for each $\e>0$, every disc $B(c,L)$ with $L>L(\e)$ contains a
point $\tau=\tau'+i\tau''$ such that the function $Q(z):\,\C\to\C$
satisfies (\ref{ap2}). Besides,  we get $|1-\exp2\pi
i\tau'|<\e,\,|1-\exp2\pi i\tau''|<\e$, hence each disc
$B(\tau,\e/8)$ intersects with $\Z^2$. Since $Q(z)$ is uniformly
continuous in $\C$, we see  that for every $\eta>0$ there exists
 a relatively dense set $E\subset\Z^2$ such that
$|Q(z+\tau)-Q(z)|<\eta$ for all $z\in\C$ and $\tau\in E$.
  Set
 $$
 q(n)=Q(n)-Q(n-1)-Q(n-i)+Q(n-1-i),\quad n\in\Z^2\subset\C.
 $$
  In view of (\ref{sp}), the function $q(n)$ has no period $\tau\in\Z^2$. Also,
take a complex number $p\in B(0,1/6)\setminus\{0\}$ and a positive
number $\g<|p|(\sup_{n\in\Z^2}|q(n)|)^{-1}$.
\medskip

{\bf Example 1}.  Put $a_n=n+p+\g q(n),\,b_n=n-p-\g
q(n),\,n\in\Z^2\subset\C$. Let us show that the divisor
$D=((a_n),\,(b_n))$ satisfies the conditions of Theorem
\ref{main}. Clearly, the divisor $D$ is almost periodic with a
regular indexing.

Condition (\ref{sep}) follows from the inequalities
$|a_n-b_{m}|>1/3$ for $n\neq m$ and $|a_n-b_n|\ge
2|p|-2\g\sup_{n\in\Z^2}|q(n)|$ for all $n$.

Furthermore, let $E$ be a convex bounded set. Since
$(a_n+b_n)/2=n$ for all $n\in\Z^2$ and $\nu_D(n)=2p+2\g q(n)$, we
get
 $$
\nu_D(E)=2p\,\card{\,E\cap\Z^2}+2\g\sum_{n\in
E}[Q(n)-Q(n-1)-Q(n-i)+Q(n-1-i)].
 $$
Clearly, $|\card{\,E\cap\Z^2}-m_2(E)|<C(m_1(\partial E)+1)$. Also,
 $$
\left|\sum_{n\in E}Q(n-1)-\sum_{n\in E}Q(n)\right|<Cm_1(\partial
E).
 $$
The same is true for the difference $\sum_{n\in
E}Q(n-i)-\sum_{n\in E}Q(n-1-i)$. Hence, we obtain (\ref{reg}) with
$\b=4\pi\bar p$.

Furthermore, for any $k\in\Z^2$ we have
 $$
\left|\sum_{1\le|n|<r}{Q(n-1-k)\over n}-
\sum_{2<|n|<r}{Q(n-k)\over n+1}\right|< 8\sup_{n\in\Z^2}|Q(n)|
 $$
 $$
  +\sup_{n\in\Z^2}|Q(n)|{\card\{n:\,|n|<r,\,|n-1|\ge r\}
+\card\{n:\,|n-1|<r,\,|n|\ge r\}\over r-2}.
 $$
Obviously, the right-hand side of this inequality is bounded
uniformly in $k\in\Z^2$ and $r>3$. By the same argument,
$$
\left|\sum_{1\le|n|<r}{Q(n-i-k)\over n}-
\sum_{2<|n|<r}{Q(n-k)\over n+i}\right|<C
 $$
and
 $$
\left|\sum_{1\le|n|<r}{Q(n-1-i-k)\over n}-
\sum_{2<|n|<r}{Q(n-k)\over n+1+i}\right|<C.
 $$
Therefore, the sum
 $$
\sum_{1\le|n|<r}{q(n-k)\over n}=\sum_{1\le|n|<r}{Q(n-k)-Q(n-1-k)-
Q(n-i-k) +Q(n-1-i-k)\over n}
 $$
 up to a bounded term has the form
 $$
 \sum_{2<|n|<r}Q(n-k)\left[{1\over n}
 -{1\over n+1}- {1\over n+i}+{1\over
 n+1+i}\right].
 $$
Note that the absolute value of the expression in the square
brackets does not exceed $C|n|^{-3}$. Taking into account the
equality $\nu^k_D(\{n\})-\nu_D(\{n\})=2\g(q(n-k)-q(n))$, we obtain
that (\ref{h2}) is satisfied for all $z\in\Z^2\subset\C$. So that
formula (\ref{f1}) with the divisor $\Dab$ and $\b=4\pi\bar p$
defines an $\AE$-function.

\medskip
{\bf Example 2}.  Put $a_n=n+p+\g q(n),\,b_n=n-p+\g q(n),
\,n\in\Z^2$. The divisor $\Dab$ is almost periodic with a regular
 indexing. Also, we obtain (\ref{sep}). Next,
$(a_n+b_n)/2=n+\g q(n),\,a_n-b_n=2p$, and for every convex bounded
set $E$ we have
$$
|\nu_D(E)-2p~\card(E\cap\Z^2)|<|\nu_D|((\partial
E)_{1/6})+2|p|~\card((\partial E)_{1/6}\cap\Z^2)<C(m_1(\partial
E)+1),
$$
where $(\partial E)_{1/6}$ is the $(1/6)$-neighborhood of $\partial
E$. So, we obtain (\ref{reg}). Finally, the difference
 $$
\sum_{1\le|n+\g q(n)|<r}{\nu_D(n+\g q(n))\over n+\g
q(n)}-\sum_{1\le|n+\g q(n-k)|<r}{\nu_D(n+\g q(n-k))\over n+\g
q(n-k)}
 $$
up to a bounded term is equal to the sum
 \begin{equation}\label{est}
\sum_{2<|n|<r}2p\g\left[{q(n-k)-q(n)\over(n+\g q(n)) (n+\g
q(n-k))}\right]
  \end{equation}
Arguing as above, we see that the sum $\sum_{2<|n|<r}q(n-k)/n^2$
is uniformly bounded  for $r>3$ and $k\in\Z^2$. Hence, the same is
true for (\ref{est}). Therefore, condition (\ref{h2}) holds for
all $z\in\Z^2$, so that formula (\ref{f1}) with $\b=4\pi\bar p$
defines  an $\AE$-function.

\medskip
{\bf Question}. If there exists an almost periodic divisor of an
$\AE$-function, which has no regular indexing?

\bigskip

School of Mathematics, Kharkov National University,

Svobody sq.,4, Kharkov 61077, Ukraine,

e-mail: Sergey.Ju.Favorov@univer.kharkov.ua
\end{document}